\documentstyle[twoside]{article}
\include{graphicx}
\title{Asymptotic and structural properties of special cases of the Wright function arising in probability theory}
\author{\sc Richard B.\ Paris \\
{\em Division of Computing and Mathematics}, \\
{\em University of Abertay Dundee, Dundee DD1 1HG, UK}\\
and\\
{\sc Vladimir Vinogradov}\\
{\em Department of Mathematics},\\
{\em Ohio University, Athens, OH 45701, USA}
\\
\\
{\em Dedicated to the memory of Lee Lorch, 1915 - 2014}
}
\begin{document}
\def\f#1#2{\mbox{${\textstyle \frac{#1}{#2}}$}}
\def\dfrac#1#2{\displaystyle{\frac{#1}{#2}}}
\def\boldal{\mbox{\boldmath $\alpha$}}
{\newcommand{\Sgoth}{S\;\!\!\!\!\!/}
\newcommand{\bee}{\begin{equation}}
\newcommand{\ee}{\end{equation}}
\newcommand{\lam}{\lambda}
\newcommand{\ka}{\kappa}
\newcommand{\al}{\alpha}
\newcommand{\fr}{\frac{1}{2}}
\newcommand{\fs}{\f{1}{2}}
\newcommand{\g}{\Gamma}
\newcommand{\br}{\biggr}
\newcommand{\bl}{\biggl}
\newcommand{\ra}{\rightarrow}
\newcommand{\mbint}{\frac{1}{2\pi i}\int_{c-\infty i}^{c+\infty i}}
\newcommand{\mbcint}{\frac{1}{2\pi i}\int_C}
\newcommand{\mboint}{\frac{1}{2\pi i}\int_{-\infty i}^{\infty i}}
\newcommand{\gtwid}{\raisebox{-.8ex}{\mbox{$\stackrel{\textstyle >}{\sim}$}}}
\newcommand{\ltwid}{\raisebox{-.8ex}{\mbox{$\stackrel{\textstyle <}{\sim}$}}}
\renewcommand{\topfraction}{0.9}
\renewcommand{\bottomfraction}{0.9}
\renewcommand{\textfraction}{0.05}
\newcommand{\mcol}{\multicolumn}
\date{}
\maketitle
\pagestyle{myheadings}
\markboth{\hfill \sc R. B.\ Paris \& V. Vinogradov  \hfill}
{\hfill \sc  Wright function asymptotics\hfill}
\begin{abstract}
This analysis paper presents previously unknown properties of some special cases of the Wright function 
whose consideration is necessitated by our work on probability theory and the theory of stochastic processes. 
Specifically, we establish new asymptotic properties of the particular Wright function
\[{}_1\Psi_1(\rho,k; \rho,0;x)= \sum_{n=0}^\infty\frac{\Gamma(k+\rho n)}{\Gamma(\rho n)}\,\frac{x^n}{n!}\qquad (|x|<\infty)\]
when the parameter $\rho\in (-1,0)\cup (0,\infty)$ and the argument $x$ is real. 
In the probability theory applications, which are focused on studies of the Poisson-Tweedie mixtures, the parameter $k$ is a non-negative integer. 
Several representations involving well-known special functions are given for certain particular values of $\rho$. The asymptotics of ${}_1\Psi_1(\rho,k;\rho,0;x)$ are obtained under numerous assumptions on the behavior of the arguments $k$ and $x$ when the parameter $\rho$ is both positive and negative.

We also provide some integral representations and structural properties involving the `reduced' Wright function
${}_0\Psi_1(-\!\!\!-; \rho,0;x)$
with $\rho\in (-1,0)\cup (0,\infty)$, which might be useful for the derivation of new properties of members of the power-variance family of distributions. Some of these imply a reflection principle that connects the functions ${}_0\Psi_1(-\!\!\!-; \pm\rho, 0;\cdot)$ and certain Bessel functions.

Several asymptotic relationships for both particular cases of this function are also given. A few of these follow under additional constraints from probability theory results which, although previously available, were unknown to analysts.
\vspace{0.4cm}

\vspace{0.3cm}

\end{abstract}

\vspace{0.3cm}

\noindent $\,$\hrulefill $\,$

\vspace{0.2cm}
\begin{center}
{\bf 1. \  Introduction}
\end{center}
\setcounter{section}{1}
\setcounter{equation}{0}
\renewcommand{\theequation}{\arabic{section}.\arabic{equation}}
In this analysis paper we deal with various asymptotic properties of a particular special function known as the Wright function that arises in analytic probability theory. This article has its outgrowth in the specialist works
of \cite{IL} and \cite{Z} and is closely related to \cite[Sections 3, 4]{PLMJ}. 

We consider the subclass of Wright functions defined by
\bee\label{e11}
{}_1\Psi_1(z)\equiv{}_1\Psi_1(\al, a;\beta, b;z)=\sum_{n=0}^\infty \frac{\g(\al n+a)}{\g(\beta n+b)}\,\frac{z^n}{n!}, 
\ee
where the parameters $\alpha$ and $\beta$ are real and positive and $a$ and $b$ are in general
arbitrary complex numbers. We also assume throughout that the $\alpha$ and $a$ are subject to the restriction
\[\alpha n+a\neq 0, -1, -2, \ldots \qquad (n=0, 1, 2, \ldots)\]
so that no gamma function in the numerator in (\ref{e11}) is singular; however, if the gamma function in the denominator is singular at the same points as the numerator, then a limiting procedure may be applied.
In the special case $\alpha=\beta=1$, the function ${}_1\Psi_1(z)$ reduces to 
the confluent hypergeometric function $(\g(a)/\g(b))\,{}_1F_1(a;b;z)$;
see, for example, \cite [ p.~40]{S}. When $\alpha=1$, the function ${}_1\Psi_1(z)$ becomes the generalized Mittag-Leffler function considered in \cite{Pr}.

We introduce the parameters associated with (\ref{e11}) given by
\bee\label{e12}
\kappa=1+\beta-\alpha, \qquad 
h=\alpha^{\alpha}\beta^{-\beta},
\qquad\vartheta=a-b.
\ee
If it is supposed that $\alpha$ and $\beta$ are such that $\kappa>0$ then ${}_1\Psi_1(z)$ 
is uniformly and absolutely convergent for all finite $z\in C$. If $\kappa=0$, the sum in (\ref{e11})
has a finite radius of convergence equal to $h^{-1}$, whereas for $\kappa<0$ the sum is divergent 
for all nonzero values of $z$. The parameter $\kappa$ is found to play a critical role 
in the asymptotic theory of ${}_1\Psi_1(z)$ by determining the sectors in the complex $z$-plane 
in which its behavior is either exponentially large or algebraic in character as $|z|\ra\infty$. 

The consideration of the entire class (\ref{e11}) of Wright functions is interesting in its own right and has already been undertaken in \cite{PLMJ, WZ1}. In contrast, in this paper we concentrate on the derivation of subtle properties of a few special cases of this function. This is partly motivated by the fact that
the function in (\ref{e11}) with $\al=\beta=\rho$, $b=0$ and $a=k$, where $k$ denotes a non-negative integer, has numerous applications  in probability theory, the 
theory of generalized linear models and 
stochastic models of actuarial science. This is because it turns out that 
the probability function of a generic member of the specific three-parameter family of univariate probability distributions 
on non-negative integers, which are  hereinafter referred to as the {\it Poisson-Tweedie mixtures}, is expressed in terms of the function ${}_1\Psi_1(z)$. We refer to 
\cite{H, JKK, BJ97, KDD, PW, PV2, PV4, V8} for more detail on this important family of discrete probability distributions. The simplest cases, which correspond to the values $\rho=-\fs$ and $\rho=1$, are considered in \cite{V8} and \cite{PV2}, respectively. 
The observation that the probability function of a generic member of this family is expressed in terms of ${}_1\Psi_1(z)$ appears to be new and will be pursued in depth in our forthcoming paper \cite{PV4}. In contrast, \cite[Section 7.5]{PW} and \cite[\S 11.1.2]{JKK} employ a double recursion for the specific probabilities associated with Poisson-Tweedie mixtures.

In this paper, we consider the particular case of the function in (\ref{e11}) given by
\bee\label{e14}
{}_1\Psi_1(\rho, k;\rho, \delta;x)=\sum_{n=0}^\infty \frac{\g(k+\rho n)}{\g(\rho n+\delta)}\,\frac{x^n}{n!}\qquad (k=0, 1, 2, \ldots\,),
\ee
where $\delta\in C$, $\rho\in (-1,0)\cup(0,\infty)$ and throughout the paper $x$ is a real variable. Although $k$ is a non-negative integer in the 
probability application, we shall occasionally relax this constraint and consider situations where $k\geq 0$ is not necessarily an integer; see, for example, (\ref{e15}).

In a similar manner to \cite[Section 2.3]{PK},  
we also define the complex-valued `reduced' Wright function $\phi(\rho, \delta, z)$ (also known as a generalized Bessel function) of argument $z \in C$,  
which is indexed by the real-valued parameter $\rho \in (-1, 0) \cup (0, \infty)$  
and the complex-valued parameter $\delta \in C$ by the following convergent series:  
\begin{equation}\label{e13}
{}_0\Psi_1(-\!\!\!-; \rho,\delta;z)\equiv
\phi(\rho, \delta; z) := ~\sum_{n = 0}^{\infty}~\frac{z^n}{n!\,\Gamma (\rho n + \delta)}~.
\end{equation}
In this paper, we shall set the parameter $\delta = 0$ in (\ref{e14}) and (\ref{e13}). In this case, the first term of the series corresponding to $n=0$ on the right-hand sides of (\ref{e14}) (when $k\neq 0$) and (\ref{e13}) is zero; compare to formula (\ref{e51}).

The derivation of the new results of the present paper, which concern numerous properties of members of the classes of the special functions (\ref{e14}) and (\ref{e13}) (with $\delta=0$), is necessitated by our studies of various approximations for some Poisson-Tweedie mixtures and members of the power-variance family of probability distributions. In view of \cite{VPY1, VPY2}, numerous representatives of this power-variance family are closely related to the `reduced' Wright function (\ref{e13}) with $\delta=0$.

At this point, we observe the following important relation for $x>0$ between the corresponding representatives of the two classes (\ref{e14}) and (\ref{e13}) of special functions introduced above. This is given by 
\begin{equation}	\label{e15}	
{}_{1}\Psi_{1}(\rho, k; \rho, 0; \epsilon x)=
\int_{0}^{\infty}~e^{-\tau}\tau^{k-1}\ \phi(\rho, 0; \epsilon x \tau^{\rho})\,d\tau,\qquad \epsilon=\mbox{sign}\,\rho
\end{equation}
valid for $k\geq 0$ (not necessarily an integer) when $\rho>0$ and for non-negative integer values of $k$ when $\rho\in (-1,0)$.
The proof of (\ref{e15}) is given in Appendix A.

The structural relationship (\ref{e15}), which connects the Wright function ${}_1\Psi_1$ with the corresponding
`reduced' Wright function $\phi$, provides the analytical explanation of the fact that the so-called Poisson-Tweedie mixtures are derived by employing the probability law of a member of the power-variance  family as the mixing measure. See also \cite[Eq.~(16)]{V8} and \cite[Eq.~(3.9)]{PV2} which pertain to the values $\rho=-\fs$ and 1, respectively.

The structure of the paper is as follows. In Section 2, some special function and polynomial representations of ${}_1\Psi_1(\rho, k;\rho, 0;x)$ are given. In Sections 3 and 4 the asymptotics of this function are considered for 
(i) $x\rightarrow\pm\infty$, $k$ finite (ii) $k\rightarrow+\infty$, $x$ finite when $\rho>0$ and $\rho\in (-1,0)$, (iii) $k\rightarrow+\infty$, $x\rightarrow+\infty$ when $\rho>0$, (iv) $k\ra+\infty$, $x\ra-\infty$ when $\rho\in(-1,0)$ and (v) $\rho\ra+\infty$ when $k$ is a positive integer and $x$ is finite. 
The asymptotics of the `reduced' Wright function $\phi(\rho, 0;x)$ for $x\ra\pm\infty$ are summarized in Section 5. The concluding Section 6 lists various representations and properties of $\phi(\rho, 0;x)$ when $\rho\in(-1,0)$. The Appendices A, B and C contain a proof of (\ref{e15}), two integral representations for ${}_1\Psi_1(\rho, k; \rho,0;x)$ stated in Section 2 as formulas (\ref{e23})--(\ref{e24}) of Theorem 1 and two specific representations of the `reduced' Wright function in terms of Bessel and Whittaker functions.

To conclude the Introduction, we point out that the majority of the results of this paper provide an analytic foundation to, or constitute extensions of the analysis counterparts of, the corresponding assertions of probability theory and the theory of stochastic processes. A few of the analytical results presented only have probabilistic proofs; in these cases, as for many of the other results presented herein, a numerical verification has been carried out.
\vspace{0.6cm}

\begin{center}
{\bf 2. \ Special representations of ${}_1\Psi_1(\rho,k;\rho,0;x)$}
\end{center}
\setcounter{section}{2}
\setcounter{equation}{0}
\renewcommand{\theequation}{\arabic{section}.\arabic{equation}}
For certain values of the parameter $\rho$ it is possible to give representations of ${}_1\Psi_1(\rho,k;\rho,0;x)$ in terms of well-known special functions. The simplest of these pertain to the values $\rho=1$ and $\rho=-\f{1}{3},\ -\fs,\ -\f{2}{3}$. See Sections 2.1--2.2 and also Sections 6.1--6.3 for the case of the `reduced' Wright function $\phi(\rho,0;x)$ when $\rho=1$ and $\rho=-\f{1}{4}$, $-\f{1}{6}$, $\f{1}{3}$, $\fs$, $\f{2}{3}$.
We remark that the representations in Sections 2.1 and 2.3 are given for real $x$, but are valid more generally for complex values of $x$.
\vspace{0.1cm}

\noindent
2.1 \ {\it Special function representations}
\vspace{0.1cm}

\noindent
When $\rho=1$ and $k>0$, we have
\begin{eqnarray}
{}_1\Psi_1(1, k; 1, 0; x)&=&\sum_{n=1}^\infty\frac{\g(k+n)\, x^n}{\g(n)\, n!}
=x\sum_{n=0}^\infty\frac{\g(k+n+1)\, x^n}{\g(n+2)\, n!}\nonumber\\
&=&x k!\,{}_1F_1(k+1;2;x),\label{e21}
\end{eqnarray}
where ${}_1F_1$ denotes the confluent hypergeometric function. In view of (\ref{e21}), numerous analytical results of our paper \cite{PV2} constitute special cases of related assertions of the present article.

Another case where a special function representation is possible corresponds to $\rho=-\fs$. In this case, we have
\begin{eqnarray}
{}_1\Psi_1(-\fs, k; -\fs, 0;\pm x)&=&\sum_{n=1}^\infty \frac{\g(k-\fs n)}{\g(-\fs n)}\,\frac{(\pm x)^n}{n!}\label{e22}\\
&=&\sum_{n=1}^\infty \frac{\g(k-n)}{\g(-n)}\,\frac{x^{2n}}{(2n)!}\pm \sum_{n=0}^\infty\frac{\g(k-n-\fs)}{\g(-n-\fs)}\,\frac{x^{2n+1}}{(2n+1)!}\equiv \Sigma_1\pm \Sigma_2.\nonumber
\end{eqnarray}
Then, for positive integer $k$, some straightforward manipulation yields
\[\Sigma_1=(-)^k \sqrt{\pi}\sum_{n=k}^\infty \frac{(\fs x)^{2n}}{\g(1-k+n) \g(n+\fs)}=(-)^k \sqrt{\pi} \sum_{m=0}^\infty \frac{(\fs x)^{2m+2k}}{m! \g(k+\fs+m)}\] 
\[=(-)^k \sqrt{\pi}\, (\fs x)^{k+1/2} I_{k-\fr}(x)\]
and 
\[\Sigma_2=\frac{-1}{\sqrt{\pi}}\sum_{n=0}^\infty\frac{(-)^n}{n!} \g(k-n-\fs) (\fs x)^{2n+1}=(-)^k \sqrt{\pi}\sum_{n=0}^\infty\frac{(\fs x)^{2n+1}}{n! \g(\f{3}{2}-k+n)}\]
\[=(-)^k \sqrt{\pi} \,(\fs x)^{k+1/2} I_{\fr-k}(x),\]
where $I_\nu(x)$ denotes the modified Bessel function of the first kind. Consequently, for non-negative integer $k$, we have the representations
\bee\label{e22a}
{}_1\Psi_1(-\fs, k;-\fs, 0;x)=(-)^k\sqrt{\pi}\,(\fs x)^{k+1/2}\{I_{k-\fr}(x)+I_{\fr-k}(x)\}
\ee
and
\bee\label{e22b}
{}_1\Psi_1(-\fs, k;-\fs, 0;-x)=\frac{2}{\sqrt{\pi}} (\fs x)^{k+1/2} K_{k-\fr}(x),
\ee
where we have employed the well-known relation $I_{-\nu}(x)=I_\nu(x)+(2/\pi) \sin (\pi\nu)\,K_\nu(x)$ connecting $I_{\pm\nu}(x)$ with the modified Bessel function of the second kind $K_\nu(x)$ (also known as the Macdonald function); see, for example, \cite[Eq.~(10.27.2)]{DLMF}.

It is worth mentioning that when $k>0$, but excluding half-integer values (for which the left-hand side of (\ref{e22}) is not defined), then $\Sigma_1\equiv 0$ and we easily find that for $x\geq 0$
\[{}_1\Psi_1(-\fs, k; -\fs, 0:\pm x)=\pm\frac{\sqrt{\pi}}{\cos \pi k} (\fs x)^{k+1/2} I_{\fr-k}(x)
\qquad (k>0;\ k\neq \fs, 1, \f{3}{2}, 2, \ldots ).\]
\vspace{0.2cm}

\noindent
2.2 \ {\it Integral representations when $\rho=-1/3$ and $\rho=-2/3$}
\vspace{0.1cm}

\noindent
The Wright function ${}_1\Psi_1(\rho,k;\rho,0;-x)$ when $\rho=-\f{1}{3}$ and $\rho=-\f{2}{3}$ is also related to special functions, but in these cases it is expressed as an infinite integral involving the $K_\nu$ Bessel and $W_{\kappa,\mu}$ Whittaker functions, respectively. 
Originally, we obtained these results by using probabilistic arguments. In this paper, we provide an analytical proof. These are given in the following theorem.
\newtheorem{theorem}{Theorem}
\begin{theorem}$\!\!\!.$\ \ For $x>0$ and positive integer $k$, we have 
\bee\label{e23}
{}_1\Psi_1(-\f{1}{3}, k; -\f{1}{3}, 0; -x)=\frac{x^{3/2}}{3\pi}\int_0^\infty e^{-t}t^{k-3/2}K_\frac{1}{3}
\left(\frac{2x^{3/2}}{3\sqrt{3t}}\right)dt
\ee
and
\bee\label{e24}
{}_1\Psi_1(-\f{2}{3}, k; -\f{2}{3}, 0; -x)=\sqrt{\frac{3}{\pi}}\int_0^\infty e^{-t}t^{k-1}\exp\left(\frac{-2x^3}{27\,t^2}\right) W_{\frac{1}{2}, \frac{1}{6}}\left(\frac{4x^3}{27\,t^2}\right)dt.
\ee
\end{theorem}

\noindent The proofs of (\ref{e23}) and (\ref{e24}) are given in Appendix B. We should point out that for $\rho=1, -\f{1}{3}, -\fs, -\f{2}{3}$, the results in (\ref{e21}), (\ref{e22b}), (\ref{e23}) and (\ref{e24}) can all be derived from (\ref{e15}); see also Section 6.1.
\vspace{0.2cm}

\noindent
2.3\ {\it Polynomial representation for $\rho\neq 0$}
\vspace{0.1cm}

\noindent
It is important that, for positive integer $k$ and 
$\rho\neq 0$, the Wright function ${}_1\Psi_1(\rho, k; \rho,0;x)$ can be expressed as
$(\rho x)^k e^x$ multiplied by a polynomial in $x^{-1}$ of degree $k-1$. This is established in (\ref{e27}), which is the main result of this subsection. 

To see this, we first note that the quotient of gamma functions 
\begin{eqnarray*}
\frac{\g(k+\rho n)}{\g(\rho n)}&=&\rho n(1+\rho n)\ldots (k-1+\rho n)=(-)^k\sum_{m=0}^k s_k^{(m)} (-\rho n)^m\\\
&=&(-)^k\sum_{m=0}^k\sum_{r=0}^k (-\rho)^m s_k^{(m)} {S}_m^{(r)} \,n(n-1)\ldots (n-r+1),
\end{eqnarray*}
where $s_k^{(m)}$ and ${S}_m^{(r)}$ are the Stirling numbers of the first and second kinds, respectively (see, for example, \cite[p.~626]{DLMF}); some Stirling numbers are given at the beginning of Section 3. 
Then
\begin{eqnarray*}
{}_1\Psi_1(\rho, k; \rho,0; x)&=&\sum_{n=1}^\infty \frac{\g(k+\rho n)}{\g(\rho n)}\,\frac{x^n}{n!}
=(-)^k\sum_{m=0}^k\sum_{r=0}^m(-\rho)^m s_k^{(m)} {S}_m^{(r)} \sum_{n=r}^\infty \frac{x^n}{(n-r)!}\\
&=& (\rho x)^k e^x \sum_{m=0}^k\sum_{r=0}^m (-\rho)^{m-k} x^{r-k} s_k^{(m)} {S}_m^{(r)},
\end{eqnarray*}
where the sum over $n$ has been evaluated as $x^re^x$.

If we put $k-r=n$ and $k-m=\ell$, where $0\leq\ell\leq k$ and $\ell\leq n\leq k$, then the above double sum can be rearranged as
\[\sum_{\ell=0}^k\sum_{n=\ell}^k (-\rho)^{-\ell} x^{-n} s_k^{(k-\ell)} {S}_{k-\ell}^{(k-n)}
=\sum_{n=0}^{k-1}x^{-n}\sum_{\ell=0}^n (-\rho)^{-\ell} s_k^{(k-\ell)} {S}_{k-\ell}^{(k-n)},\]
where, since $s_m^{(0)}={S}_m^{(0)}=\delta_{0,m}$, the inner sum vanishes when $n=k$.  
If we define the polynomial ${\bf h}_m(u)$ of degree $m$ by
\bee\label{e27a}
{\bf h}_m(u):=\sum_{n=0}^m {D}_n u^n,\qquad 
{D}_n:=\sum_{\ell=0}^n (-)^{\ell}\rho^{n-\ell} s_k^{(k-\ell)} {S}_{k-\ell}^{(k-n)},
\ee
then, for all $x$ (both positive and negative) and $\rho\neq 0$, we can write
\bee\label{e27}
{}_1\Psi_1(\rho, k; \rho,0; x)=(\rho x)^k e^x {\bf h}_{k-1}((\rho x)^{-1}).
\ee 
See also (\ref{e33}) for the first few coefficients $D_n$. The representation (\ref{e27}) generalizes \cite[Eq.~(2.5)]{PV2} where the case $\rho=1$ was considered.

To conclude, we list the specific examples corresponding to $1\leq k\leq 5$:
\begin{eqnarray*}
{}_1\Psi_1(\rho, 1; \rho,0;x)&=&\rho xe^x,\\
{}_1\Psi_1(\rho, 2; \rho,0;x)&=&(\rho x)^2e^x\left\{1+\frac{1+\rho}{\rho x}\right\},\\
{}_1\Psi_1(\rho, 3; \rho,0;x)&=&(\rho x)^3e^x\left\{1+\frac{3(1+\rho)}{\rho x} +\frac{(1+\rho)(2+\rho)}{(\rho x)^2}\right\},\\
{}_1\Psi_1(\rho, 4; \rho,0;x)&=&(\rho x)^4 e^x\left\{1+\frac{6(1+\rho)}{\rho x}+\frac{(1+\rho)(11+7\rho)}{(\rho x)^2}+\frac{(1+\rho)(2+\rho)(3+\rho)}{(\rho x)^3}\right\},\\
{}_1\Psi_1(\rho, 5; \rho,0;x)&=&(\rho x)^5 e^x\left\{1+\frac{10(1+\rho)}{\rho x}+\frac{5(1+\rho)(7+5\rho)}{(\rho x)^2}+\frac{5(1+\rho)(2+\rho)(5+3\rho)}{(\rho x)^3}\right.\\
&&\hspace{6cm}\left. +\frac{(1+\rho)(2+\rho)(3+\rho)(4+\rho)}{(\rho x)^4}\right\}.
\end{eqnarray*}

\vspace{0.6cm}

\begin{center}
{\bf 3. \ Asymptotics of ${}_1\Psi_1(\rho,k;\rho,0;x)$ for fixed $k$}
\end{center}
\setcounter{section}{3}
\setcounter{equation}{0}
\renewcommand{\theequation}{\arabic{section}.\arabic{equation}}
\vspace{0.2cm}

\noindent
3.1 \ {\it Asymptotics for $x\ra\pm\infty$ when $k$ is a positive integer}
\vspace{0.1cm}

\noindent
When $k$ is a fixed positive integer and $\rho\neq 0$, the asymptotic expansion of ${}_1\Psi_1(\rho,k;\rho,0;x)$ as $x\ra\pm\infty$ follows immediately from (\ref{e27}) as
\bee\label{e31}
{}_1\Psi_1(\rho,k;\rho,0;x)=(\rho x)^k e^x \sum_{n=0}^{k-1} {D}_n (\rho x)^{-n}\qquad (k\in {\bf N}),
\ee
where the coefficients ${D}_n$ are defined in (\ref{e27a}).

To determine the first few coefficients ${D}_n$ we note the following evaluations of the Stirling numbers of the first and second kinds:
\[s_n^{(n)}=1,\quad s_n^{(n-1)}=-\left(\!\!\begin{array}{c}n\\2\end{array}\!\!\right),\quad
s_n^{(n-2)}=\frac{1}{4}(3n-1)\left(\!\!\begin{array}{c}n\\3\end{array}\!\!\right),\quad
s_n^{(n-3)}=-\left(\!\!\begin{array}{c}n\\2\end{array}\!\!\right)\left(\!\!\begin{array}{c}n\\4\end{array}\!\!\right),\]
\[{S}_n^{(n)}=1,\quad {S}_n^{(n-1)}=\left(\!\!\begin{array}{c}n\\2\end{array}\!\!\right),\quad
{S}_n^{(n-2)}=\frac{1}{4}(3n-5)\left(\!\!\begin{array}{c}n\\3\end{array}\!\!\right),\quad
{S}_n^{(n-3)}=\left(\!\!\begin{array}{c}n-2\\2\end{array}\!\!\right)\left(\!\!\begin{array}{c}n\\4\end{array}\!\!\right).\]
Then some routine algebra yields that
\[{D}_0=1,\quad {D}_1=(1+\rho)\left(\!\!\begin{array}{c}k\\2\end{array}\!\!\right),\quad
{D}_2=\frac{1}{4}(1+\rho) \left(\!\!\begin{array}{c}k\\3\end{array}\!\!\right)\{3k(1+\rho)-1-5\rho\},\]
\bee\label{e33}
{D}_3=\frac{1}{2}(1+\rho) \left(\!\!\begin{array}{c}k\\4\end{array}\!\!\right) \{k(1+\rho)+2\rho\}\{k(1+\rho)+1+3\rho\}, \ldots\ .\ee

We make three remarks concerning formula (\ref{e31}). First, it is an {\it exact\/} result. Secondly, this expansion consists of a finite series in inverse powers of $x$ valid for a fixed admissible $\rho\neq 0$ (that is, positive or negative $\rho$). And thirdly, it holds for both $x\ra+\infty$ and $x\ra-\infty$. Consequently, it follows that ${}_1\Psi_1(\rho,k;\rho,0;x)$ grows exponentially for $x>0$ and decays exponentially for $x<0$. 
\vspace{0.2cm}

\noindent
3.2 \ {\it Asymptotics for $\rho\ra+\infty$ when $k$ is a positive integer}
\vspace{0.1cm}

\noindent
An approximation for ${}_1\Psi_1(\rho, k; \rho,0;x)$ for $\rho\ra+\infty$ when $k$ is a fixed positive integer
and $x$ is fixed (positive or negative) also follows directly from (\ref{e27a}) and (\ref{e27}). We have
\begin{eqnarray}
{}_1\Psi_1(\rho, k;\rho, 0;x)&=&\rho^k e^x \sum_{n=0}^{k-1} x^{k-n} \sum_{\ell=0}^n \frac{(-)^\ell}{\rho^\ell}\,
s_k^{(k-\ell)} {S}_{k-\ell}^{(k-n)}\nonumber\\
&\sim& \rho^k e^x \sum_{n=0}^{k-1} x^{k-n} \bl\{S^{(k-n)}_k + \frac{k(k-1)}{2\rho}\,S^{(k-n)}_{k-1}+
O(\rho^{-2})\br\} \label{e33b}
\end{eqnarray}
as $\rho\ra+\infty$. More terms can be easily generated if required.

\vspace{0.2cm}

\noindent
3.3 \ {\it Asymptotics for $x\ra\pm\infty$ when $k>0$}
\vspace{0.1cm}

\noindent
When the restriction on $k$ being a positive integer is relaxed to $k>0$, the expansion in (\ref{e31}) ceases to be valid. In this case, we can apply the asymptotic theory of integral functions of hypergeometric type; see, for example,
\cite{Br}, \cite[\S 2.3]{PK}. For the Wright function in (\ref{e11}) with the numerator parameters $\alpha$, $a$ and denominator parameters $\beta$, $b$, we have 
the following (formal) exponential and algebraic asymptotic expansions $E_{1,1}(z)$ and $H_{1,1}(z)$ given by
\bee\label{e20a}
E_{1,1}(z):=Z^\vartheta e^Z\sum_{j=0}^\infty A_jZ^{-j},\qquad Z=\kappa(hz)^{1/\kappa}
\ee
and 
\bee\label{e20b}
H_{1,1}(z):=\frac{1}{\alpha}\sum_{s=0}^\infty \frac{(-)^s}{s!}\,\frac{\g((s+a)/\alpha)}{\g(b-\beta(s+a)/\alpha)}\,
z^{-(s+a)/\alpha},
\ee
where the quantities $\kappa$, $h$ and $\vartheta$ are defined in (\ref{e12}); see \cite[pp.~56, 57]{PK}.
The coefficients $A_j$ in the exponential expansion (\ref{e20a}) appear in the inverse factorial expansion of the ratio of gamma functions 
$\g(1-b+\beta s)/\g(1-a+\alpha s)$ for large $|s|$,
with 
\[A_0=\kappa^{-\vartheta-1/2}\alpha^{a-1/2}\beta^{1/2-b};\]
an algorithm for their determination is described in \cite[pp.~46--48]{PK}.
Then, if $0<\kappa<2$, we have
\bee\label{e34}
{}_1\Psi_1(\alpha,a; \beta;b; z)\sim\left\{\begin{array}{ll} E_{1,1}(z)+H_{1,1}(ze^{\mp\pi i}) & (|\arg\,z|\leq\fs\pi\kappa)\\
\\
H_{1,1}(ze^{\mp\pi i}) & (|\arg (-z)|<\fs\pi(2-\kappa))\end{array}\right.
\ee
as $|z|\ra\infty$, where the upper or lower sign in $H_{1,1}(ze^{\mp\pi i})$ is chosen according as $z$ lies in the upper or lower half-plane, respectively. 
A more precise result that takes into account the Stokes phenomenon on the rays $\arg\,z=\pm\pi\kappa$ when $0<\kappa\leq 1$ is given in \cite{PLMJ}.

We now specialize the result in (\ref{e34}) to derive the asymptotics of ${}_1\Psi_1(\rho,k; \rho,0; x)$ for large $|x|$ in the case where $\rho>0$ and $k>0$. We have
\[{}_1\Psi_1(\rho,k; \rho,0; x)=\sum_{n=0}^\infty \frac{\g(k+\rho n)}{\g(\rho n)}\,\frac{x^n}{n!},\]
so that
$\kappa=h=1$, $\vartheta=k$, $A_0=\rho^k$. Then we obtain\footnote{The subdominant algebraic expansion undergoes a Stokes phenomenon on the positive real $x$-axis and is omitted.} from (\ref{e20a})--(\ref{e34}) that
\bee\label{e35}
{}_1\Psi_1(\rho,k; \rho,0; x)\sim\left\{\begin{array}{ll}\displaystyle{x^ke^x\sum_{j=0}^\infty  A_jx^{-j}} & (x\ra+\infty)\\
\displaystyle{\frac{1}{\rho}\sum_{s=0}^\infty \frac{(-)^s}{s!}\,\frac{\g((s+k)/\rho)}{\g(-s-k)}
(-x)^{-(s+k)/\rho}+x^ke^x\sum_{j=0}^\infty  A_jx^{-j}} & (x\ra-\infty) \end{array}\right.
\ee
for fixed $k>0$, where the coefficients $A_j=\rho^{k-j} {D}_j$ with the $D_j$ defined in (\ref{e27a}).

The expansion of ${}_1\Psi_1(\rho,k; \rho,0; x)$ is seen to be dominated by the algebraic expansion as $x\ra-\infty$. 
We observe that this last expansion vanishes when $k$ is a positive integer and note the fact that, in this case, the coefficients $A_j=0$ for $j\geq k$. The expansion in (\ref{e35}) as $x\ra\pm\infty$ then agrees with that in (\ref{e31}).

\vspace{0.6cm}

\begin{center}
{\bf 4. \ Asymptotics of ${}_1\Psi_1(\rho,k;\rho,0;x)$ for $k\ra+\infty$}
\end{center}
\setcounter{section}{4}
\setcounter{equation}{0}
\renewcommand{\theequation}{\arabic{section}.\arabic{equation}}
In this section, we consider the asymptotics of ${}_1\Psi_1(\rho,k; \rho,0;x)$ for $k\ra+\infty$ when $x$ is finite (or such that $x=o(k)$) and also when $x$ is allowed to become large, namely of order $k$. The cases $\rho>0$ and $\rho\in (-1,0)$ are considered separately as their treatments are distinct. Apart from the situation corresponding to $\rho\in (-1,0)$ with $x$ finite, the approach we employ is essentially the saddle-point method. This is consistent with the fact that the case of negative values of $\rho$ corresponds to the second polar type of the mechanism of formation of the probabilities of large deviations, whereas it is their Cram\'er's or first type which is regarded as the probabilistic counterpart of the saddle-point method; see \cite[Introduction]{V94} for more details.
\vspace{0.4cm}

\noindent
4.1\ {\it The case $\rho>0$}
\vspace{0.1cm}

\noindent
An integral representation suitable for the determination of the large-$k$ asymptotics of the function ${}_1\Psi_1(\rho,k; \rho,0;x)$ when $\rho>0$ can be obtained as follows.
From the loop integral representation for the beta function \cite[Eq.~(5.12.10)]{DLMF}, we see that
\[\frac{\g(k+\rho n)}{k! \g(\rho n)}=\frac{1}{2\pi i} \int_0^{(1+)}\!\!\!\frac{t^{k+\rho n-1}}{(t-1)^{k+1}}dt\qquad (k+\rho n>0),\]
where $k>0$ (not necessarily an integer) and the integration path is a loop that starts at the origin, encircles the point $t=1$ in the positive sense and returns to the origin. Then, substitution of this result into (\ref{e14}) followed by reversal of the order of summation and integration yields 
the integral representation
\begin{eqnarray}
{}_1\Psi_1(\rho,k;\rho,0;x)&=&\sum_{n=1}^\infty \frac{\g(k+\rho n)}{\g(\rho n)}\,\frac{x^n}{n!}
=\frac{k!}{2\pi i}\int_0^{(1+)}\!\!\!\frac{t^{k-1}}{(t-1)^{k+1}}\,(e^{xt^\rho}-1)\,dt\nonumber\\
&=&\frac{k!}{2\pi i}\int_0^{(1+)}\!\!\!\frac{t^{k-1}\,e^{xt^\rho}}{(t-1)^{k+1}}\,dt,\label{e400}
\end{eqnarray}
since $t^{k-1}/(t-1)^{k+1}$ has zero residue at $t=1$.
\vspace{0.2cm}

\noindent
4.1.1\ {\it Asymptotics for $k\ra+\infty$ and $x\ra+\infty$}
\vspace{0.1cm}

\noindent
To enable an estimation of ${}_1\Psi_1(\rho,k; \rho,0;x)$ for $k\ra+\infty$ with $x=O(k)$, we put $x=ku$, where $u>0$ is fixed. In view of (\ref{e400}), we have
\bee\label{e43}
{}_1\Psi_1(\rho,k;\rho,0;x)=\frac{k!}{2\pi i}\int_0^{(1+)}\!\!\!\frac{e^{k\psi(t)}}{t(t-1)}\,dt,\qquad \psi(t)=\log\,\frac{t}{t-1}+ut^\rho.
\ee

The function $\psi(t)$ has saddle points $t_s$ where $\psi'(t)=0$; that is when
\bee\label{e44}
t_s^\rho(t_s-1)=\frac{1}{\rho u}.
\ee
When $0<\rho<1$, there is a single saddle\footnote{It is worth mentioning that an application of \cite[Proposition 5]{BKKK} implies (with some effort) that given $\rho > 0$, the solution $t_{s0}(u)$ to 
(\ref{e44})  
admits the following closed-form representation in terms of the `reduced' Wright function $\phi$:
\[t_{s0}(u) = \bl\{ u 
- \log\bl(~ \int_0^{\infty} \frac{e^{-u y}}{y (1 + y)} \, 
\phi \bl(-\frac{\rho}{\rho + 1}, 0; - \frac{1 + y}{y ^{\rho/(\rho + 1)}} \br) dy \br)
\br\}^{1/(\rho + 1)}.
\]} 
at $t_{s0}$ on the positive real $t$-axis satisfying $t_{s0}>1$. The steepest descent path through this saddle crosses the real axis in a perpendicular direction and starts and finishes at the origin; see Fig.~1(a). In this case, the integration path can be made to coincide with the steepest descent path through $t_{s0}$. When $\rho=1$, an additional saddle moves off an adjacent Riemann sheet and appears on the negative $t$-axis. 
When $\rho>1$, additional saddles arise, which are situated symmetrically about the real $t$-axis on an
approximately circular path surrounding the interval $[0,1]$. The directions of the valleys of $\psi(t)$ (where $\psi(t)\ra-\infty$ as $|t|\ra\infty$) are easily seen to be the rays $\arg\,t=\pm(2s+1)\pi/\rho$, $s=0, 1, \ldots\, $;
the directions of the hills (or ridges) are symmetrically positioned between consecutive valleys on the rays $\arg\,t=\pm 2\pi s/\rho$. Paths of steepest descent can only pass to infinity down one of the valleys or terminate at $t=0$, whereas steepest ascent paths pass to infinity on a ridge or terminate at the singular point $t=1$.

A study of the steepest descent paths through these saddles reveals that the integration path in (\ref{e43}) can be deformed to pass over those {\it three} saddles with the greatest real parts as shown in Fig.~1. The dominant saddle $t_{s0}>1$ is on the positive real $t$-axis, with the two neighboring complex saddles $t_{s,\pm 1}$ being subdominant, since it is reasonably easy to show that $\Re (\psi(t_{s,\pm 1}))<\psi(t_{s0})$ (we omit these details).
\begin{figure}[th]
	\begin{center}	{\tiny($a$)}\includegraphics[width=0.3\textwidth]{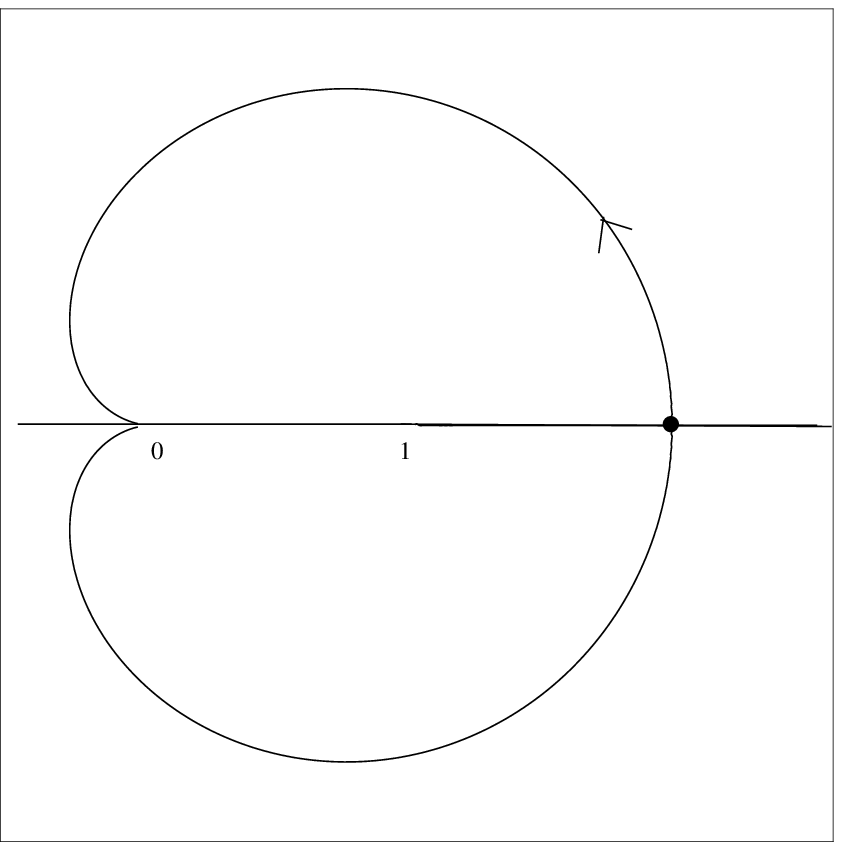}\qquad
	{\tiny($b$)}\includegraphics[width=0.3\textwidth]{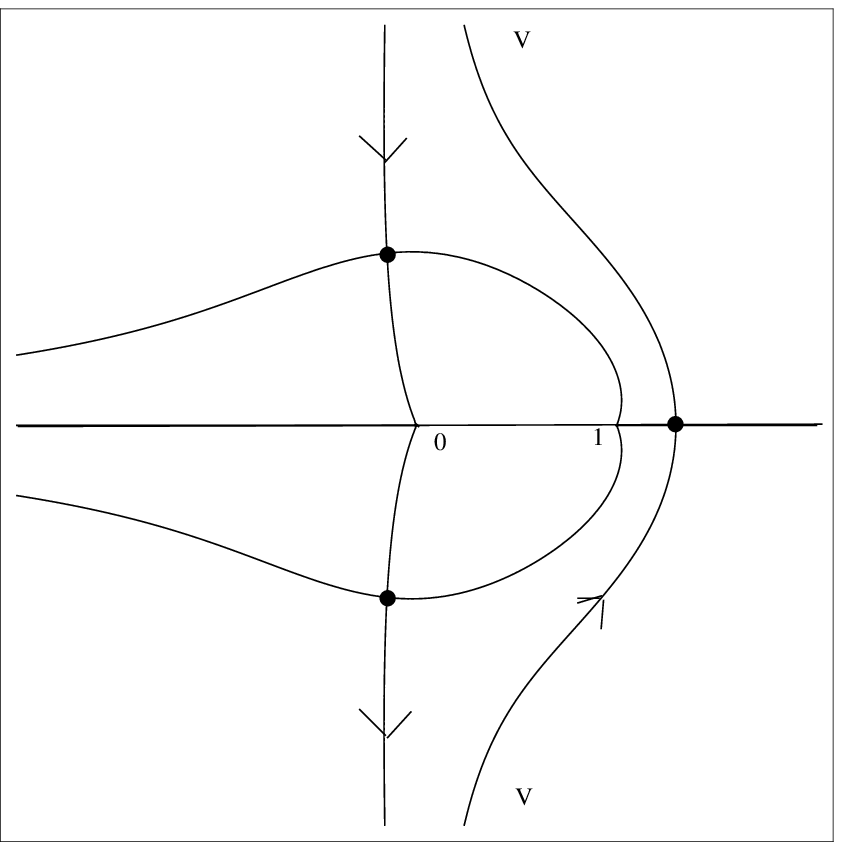}\\
	\vspace{0.3cm}
	
	{\tiny($c$)}\includegraphics[width=0.3\textwidth]{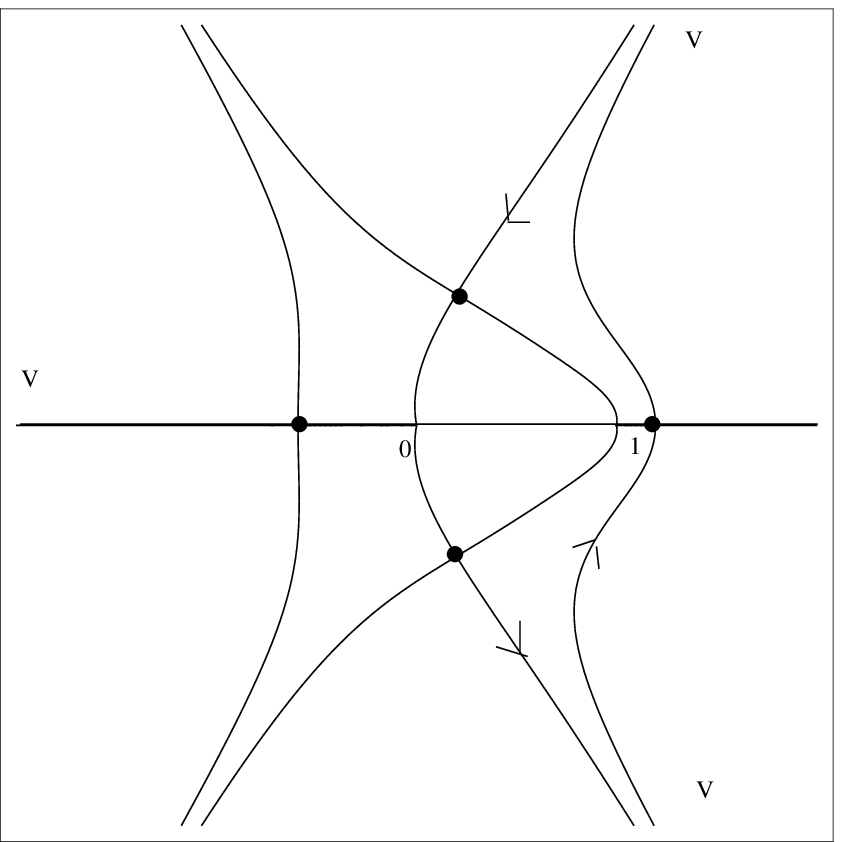}\qquad
	{\tiny($d$)}\includegraphics[width=0.3\textwidth]{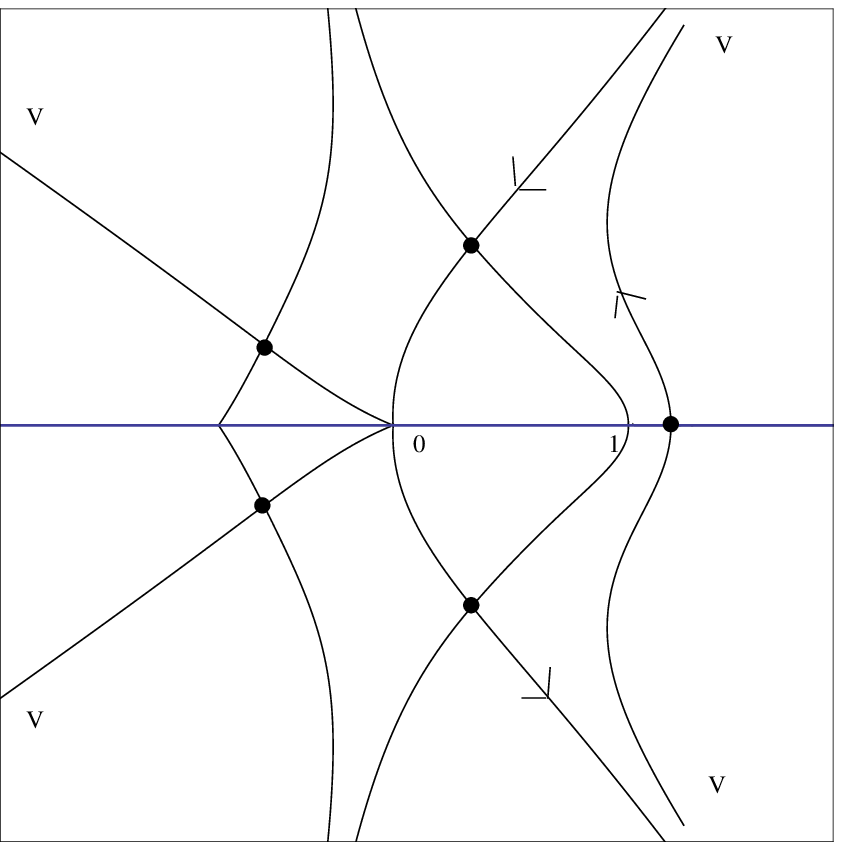}
\caption{\small{Typical paths of steepest descent and ascent through the saddles when (a) $\rho=2/3$, (b) $\rho=2$, (c) $\rho=3$ and (d) $\rho=3.6$. The saddles are denoted by heavy dots and the valleys by V; the arrows indicate the direction of integration in (\ref{e43}) taken along steepest descent paths. There is a branch cut along $(-\infty,1]$.}}
	\end{center}
\end{figure}

Straightforward differentiation shows that at a saddle $t_s$,
\[\psi''(t_s)=\frac{2t_s-1}{t_s^2(t_s-1)^2}+\rho(\rho-1)\,ut_s^{\rho-2}
=\frac{(1+\rho)}{t_s^2(t_s-1)^2}\,\left(t_s-\frac{\rho}{1+\rho}\right)\]
upon use of (\ref{e44}) to simplify the second term. Then a routine application of the saddle-point method applied to the dominant saddle $t_{s0}$  (we ignore the subdominant contributions present when $\rho\geq 1$) produces 
\[{}_1\Psi_1(\rho,k;\rho,0;ku)=\frac{k!}{2\pi} \sqrt{\frac{2\pi}{k\psi''(t_{s0})}}\,\frac{\exp\,[k\psi(t_{s0})]}{t_{s0}(t_{s0}-1)}\left\{1+O\left(k^{-1}\right)\right\}\]
\bee\label{e45}
=\g(k)\sqrt{\frac{k}{2\pi(1+\rho)}} \left(t_{s0}-\frac{\rho}{1+\rho}\right)^{\!\!-1/2} \!\exp\, [k\psi(t_{s0})]\,\{1+O(k^{-1})\}
\ee
as $k\ra\infty$, where from (\ref{e43})
\bee\label{e46}
\psi(t_{s0})=-\log\,(1-\tau_s)+\frac{\tau_s}{\rho(1-\tau_s)},\qquad \tau_s=1/t_{s0}.
\ee
The value of the dominant saddle $t_{s0}\equiv t_{s0}(u)>1$ is determined from (\ref{e44}). 
It is clear from (\ref{e46}) that $\psi(t_{s0})>0$, with the result that the factor $\exp\,[k\psi(t_{s0})]$ is exponentially large as $k\ra+\infty$.
 
In Table 1 we present some numerical values for the absolute relative error in the computation of ${}_1{\hat \Psi}_1\equiv {}_1\Psi_1(\rho,k;\rho,0;ku)/\g(k)$ using the leading asymptotic behavior in (\ref{e45}) for different $\rho$ and $u$ and varying values of $k$. 
\begin{table}[tb]
\caption{\footnotesize{Values of the absolute relative error in the computation of ${}_1{\hat \Psi}_1$ for different values of $k$ obtained from (\ref{e45}). }}
\begin{center}
\begin{tabular}{r|ll|ll}
\hline
&&&\\[-0.25cm]
\mcol{1}{c|}{} & \mcol{2}{c|}{$u=1,\ \ \rho=\fs$} & \mcol{2}{c}{$u=\fs,\ \ \rho=2$}  \\ 
\mcol{1}{c|}{$k$} & \mcol{1}{c}{${}_1{\hat \Psi}_1$} & \mcol{1}{c|}{Error}  & \mcol{1}{c}{${}_1{\hat \Psi}_1$} & \mcol{1}{c}{Error}\\
[.1cm]\hline
&&&\\[-0.3cm]
20  & $1.373292\times 10^{18}$ & $1.237\times 10^{-2}$ & $2.215937\times 10^{19}$ & $1.084\times 10^{-2}$ \\
30  & $1.957497\times 10^{27}$ & $8.192\times 10^{-3}$ & $1.200853\times 10^{29}$ & $7.180\times 10^{-3}$ \\
50  & $8.405994\times 10^{45}$ & $4.889\times 10^{-3}$ & $3.021946\times 10^{48}$ & $4.286\times 10^{-3}$ \\
80  & $6.720337\times 10^{72}$ & $3.047\times 10^{-3}$ & $3.280812\times 10^{77}$ & $2.671\times 10^{-3}$ \\
100 & $1.009961\times 10^{91}$ & $3.435\times 10^{-3}$ & $7.133294\times 10^{96}$ & $2.135\times 10^{-3}$ \\
[.2cm]\hline
\end{tabular}
\end{center}
\end{table}
\vspace{0.2cm}

\noindent
4.1.2\ {\it Asymptotics for $k\ra+\infty$ and $k/x\ra+\infty$}
\vspace{0.1cm}

\noindent
When $k\ra+\infty$ and $x=o(k)$ (which includes the case of finite $x=O(1)$), 
we find from (\ref{e43}) or (\ref{e44}) the equation for the saddle points in the form
\[t_s^\rho(t_s-1)=\frac{k}{\rho x}.\]
The dominant saddle on the positive real $t$-axis consequently has the approximate value 
\[t_{s0}\simeq (k/(\rho x))^{1/(1+\rho)}\qquad (k\ra+\infty).\]
More precisely, a routine perturbation expansion procedure shows that
\bee\label{e47}
t_{s0}=\chi \sum_{r=0}^\infty B_r \chi^{-r},\qquad \chi:=\left(\frac{k}{\rho x}\right)^{\!1/(1+\rho)},
\ee
where
\[B_0=1,\quad B_1=\frac{1}{1+\rho},\quad B_2=\frac{\rho}{2(1+\rho)^2},\quad B_3=\frac{\rho(\rho-1)}{3(1+\rho)^3},\]
\[B_4=\frac{\rho(2-5\rho+2\rho^2)}{8(1+\rho)^4}, \quad B_5=\frac{\rho}{15(1+\rho)^3} (-3+13\rho-13\rho^2+3\rho^3),\]
\[B_6=\frac{\rho}{144(1+\rho)^6} (24-154\rho+269\rho^2-154\rho^3+24\rho^4), \ldots\ .\]

Substitution of (\ref{e47}) into (\ref{e46}) expressed in the form 
\[\psi(t_{s0})=\sum_{r=1}^\infty t_{s0}^{-r} \left(\frac{1}{r}+\frac{1}{\rho}\right)\qquad (t_{s0}>1)\]
yields, with the help of {\it Mathematica}, that
\[\psi(t_{s0})
=\frac{(1+\rho)}{\rho}\chi^{-1}+\frac{1}{2}\chi^{-2}+\frac{2\rho-1}{6(1+\rho)} \chi^{-3}
+\frac{(3\rho-1)(\rho-1)}{12(1+\rho)^2} \chi^{-4}\]
\[+\frac{(4\rho-1)(3\rho-2)(2\rho-3)}{120(1+\rho)^3} \chi^{-5}
+\frac{(5\rho-1)(2\rho-1)(\rho-1)(\rho-2)}{60(1+\rho)^4}\chi^{-6}+O(\chi^{-7})\]
for $k/x\ra+\infty$.
From (\ref{e45}), we then obtain the leading behavior
\bee\label{e48}
{}_1\Psi_1(\rho,k;\rho,0;x)=\frac{\g(k)}{\sqrt{2\pi(1+\rho)}} (\rho xk^\rho)^{1/(2(1+\rho))}\, \exp [k\psi(t_{s0})]
\,\{1+O(k^{-\rho/(1+\rho)})\},
\ee
as $k\ra+\infty$ when $k/x\ra+\infty$, where
\bee\label{e48a}
k\psi(t_{s0})=
(\rho k^\rho x)^{1/(1+\rho)}\bl\{\frac{1+\rho}{\rho}+\frac{1}{2} \chi^{-1} 
+\frac{(2\rho-1)}{6(1+\rho)} \chi^{-2}
+\frac{(3\rho-1)(\rho-1)}{12(1+\rho)^2} \chi^{-3}+O(\chi^{-4})\br\}.
\ee

\begin{table}[t]
\caption{\footnotesize{Values of the absolute relative error in the computation of ${}_1{\hat \Psi}_1$ for different values of $k$ obtained from (\ref{e48}) and (\ref{e48a}). }}
\begin{center}
\begin{tabular}{r|ll|ll}
\hline
&&&\\[-0.25cm]
\mcol{1}{c|}{} & \mcol{2}{c|}{$x=1,\ \ \rho=\fs$} & \mcol{2}{c}{$x=1,\ \ \rho=\f{3}{2}$}  \\ 
\mcol{1}{c|}{$k$} & \mcol{1}{c}{${}_1{\hat \Psi}_1$} & \mcol{1}{c|}{Error}  & \mcol{1}{c}{${}_1{\hat \Psi}_1$} & \mcol{1}{c}{Error}\\
[.1cm]\hline
&&&\\[-0.3cm]
20   & $6.966593\times 10^{01}$ & $1.117\times 10^{-1}$ & $3.665626\times 10^{05}$ & $3.207\times 10^{-3}$ \\
50   & $5.142250\times 10^{02}$ & $7.559\times 10^{-2}$ & $3.555739\times 10^{09}$ & $6.793\times 10^{-3}$ \\
100  & $3.547863\times 10^{03}$ & $5.714\times 10^{-2}$ & $2.108775\times 10^{14}$ & $7.198\times 10^{-3}$ \\
200  & $3.909080\times 10^{04}$ & $4.364\times 10^{-2}$ & $3.015835\times 10^{21}$ & $6.738\times 10^{-3}$ \\
500  & $2.368065\times 10^{06}$ & $3.093\times 10^{-2}$ & $5.675040\times 10^{36}$ & $5.625\times 10^{-3}$ \\
[.2cm]\hline
&&&\\[-0.25cm]
\mcol{1}{c|}{} & \mcol{2}{c|}{$x=10,\ \ \rho=\fs$} & \mcol{2}{c}{$x=10,\ \ \rho=\f{3}{2}$}  \\ 
\mcol{1}{c|}{$k$} & \mcol{1}{c}{${}_1{\hat \Psi}_1$} & \mcol{1}{c|}{Error}  & \mcol{1}{c}{${}_1{\hat \Psi}_1$} & \mcol{1}{c}{Error}\\
[.1cm]\hline
&&&\\[-0.3cm]
50   & $3.511690\times 10^{14}$ & $5.279\times 10^{-2}$ & $2.620759\times 10^{27}$ & $1.794\times 10^{-1}$ \\
100  & $1.403708\times 10^{18}$ & $3.454\times 10^{-2}$ & $3.486374\times 10^{39}$ & $1.307\times 10^{-1}$ \\
200  & $5.198767\times 10^{22}$ & $2.311\times 10^{-2}$ & $5.286625\times 10^{57}$ & $9.396\times 10^{-2}$ \\
500  & $4.579438\times 10^{30}$ & $1.399\times 10^{-2}$ & $3.054764\times 10^{96}$ & $6.012\times 10^{-2}$ \\
1000 & $3.340732\times 10^{38}$ & $9.751\times 10^{-3}$ & $2.744393\times 10^{143}$& $4.273\times 10^{-2}$ \\
[.2cm]\hline
\end{tabular}
\end{center}
\end{table}

In the cases $\rho=1$ (where we note that the term involving $\chi^{-3}$ in (\ref{e48a}) vanishes) and $\rho=2$, we consequently find from (\ref{e47})--(\ref{e48a}) the approximations
\[{}_1\Psi_1(1,k;1,0;x)= \frac{\g(k)}{2\sqrt{\pi}}\,(kx)^{1/4}\hspace{8.5cm} \]
\bee\label{e49}
\hspace{1cm}\times\exp\,[2(kx)^{1/2}\{1+\f{1}{4}(x/k)^{1/2}+\f{1}{24}(x/k)+O((x/k)^2)\}]\,\{1+O(k^{-1/2})\}
\ee
and
\[{}_1\Psi_1(2,k;2,0;x)= \frac{\g(k)}{\sqrt{6\pi}}\,(2k^2x)^{1/6} \hspace{8.5cm}\]
\bee\label{e410}
\hspace{1cm}\times\exp\,[\f{3}{2}k^{2/3} (2x)^{1/3}\{1+\f{1}{3}(2x/k)^{1/3}+\f{1}{9}(2x/k)^{2/3}+O(x/k)\}]\,\{1+O(k^{-2/3})\}
\ee
as $k\ra+\infty$ and $k/x\ra+\infty$. 

From (\ref{e21}), the case $\rho=1$ can also be expressed in terms of the confluent hypergeometric function ${}_1F_1(1+k;2;x)$. The asymptotic behavior of this function follows with some effort from \cite[Eq.~(3.8.3)]{S} as
\[{}_1F_1(1+k;2;x) = \frac{(kx)^{-3/4}}{2\sqrt{\pi}}\,\exp\,[2(kx)^{1/2}+\fs x]\,\{1+ \f{1}{12}(x^3/k)^{1/2}+O(x^3/k)\}\]
as $k\ra+\infty$ with $k/x\ra+\infty$ (see also \cite[Theorem 2.3]{PV2} where a more accurate asymptotic expansion than this formula is stated for fixed $x>0$). Combination of the above result with (\ref{e21}) then yields (\ref{e49}), albeit with a different order term. A related result on the double asymptotics of the function ${}_1F_1$ can also be found in 
\cite[Theorem~2.3.i]{PV2}. Note that \cite[Corollary~2.4]{PV2} addresses the case where $\rho=1$, $k\ra+\infty$ and $x=O(k)$, which is not considered in the present paper.

In Table 2 we show values of the absolute relative error in the computation of the normalized function ${}_1{\hat \Psi}_1\equiv {}_1\Psi_1(\rho,k;\rho,0;x)/\g(k)$ using the leading asymptotic behavior in (\ref{e48}) and (\ref{e48a}) for different $\rho$ and $x$ and varying values of $k$. 

\vspace{0.4cm}

\noindent
4.2\ {\it The case $\rho\in (-1,0)$}
\vspace{0.1cm}

\noindent
Throughout this subsection we consider $\rho\in (-1,0)$ and  write $\rho=-\sigma$ with $0<\sigma<1$. From (\ref{e31}) and (\ref{e33}),
the behavior of ${}_1\Psi_1(-\sigma, k; \sigma, 0;-x)$  is given by
\[{}_1\Psi_1(-\sigma, k; \sigma, 0;-x)\sim (\sigma x)^k e^{-x}\bl\{1+\frac{(1-\sigma)k(k-1)}{2\sigma x}+\hspace{5cm}\]
\bee\label{e420}
\hspace{5cm}+\frac{(1-\sigma)k(k-1)(k-2)}{24(\sigma x)^2}\,\{3k(1-\sigma)+5\sigma-1\}+\ldots\br\}
\ee
for $x\ra+\infty$ and fixed positive integer $k$.
From the structure of the expansion in (\ref{e31}) and the closed-form expressions (\ref{e33}) for the initial coefficients $D_n$, it follows that (\ref{e420}) with the first three terms in curly braces is {\it exact\/} for $k=1, 2$ and $3$.
\vspace{0.2cm}

\noindent
4.2.1\ {\it Asymptotics for $k\ra+\infty$ and finite $x$}
\vspace{0.1cm}

\noindent
We have for $k>0$ (not necessarily an integer) and $\rho=-\sigma$ that
\[{}_1\Psi_1(\rho,k;\rho,0;x)=\sum_{n=1}^\infty \frac{\g(k-\sigma n)}{\g(\rho n)}\,\frac{x^n}{n!}\]
\[=\frac{x\g(k)}{\g(\rho)}\left\{\frac{\g(k-\sigma)}{\g(k)}+\frac{x\g(\rho)}{2!\g(2\rho)}\,
\frac{\g(k-2\sigma)}{\g(k)}+
\frac{x^2\g(\rho)}{3!\g(3\rho)}\,\frac{\g(k-3\sigma)}{\g(k)}+ \cdots \right\}.\]
As $k\ra+\infty$, the ratio $\g(k-m\sigma)/\g(k-(m+1)\sigma)\sim k^\sigma$ for finite positive integer $m$, with the consequence that
\[\g(k-\sigma)\gg \g(k-2\sigma) \gg \g(k-3\sigma) \ldots\ .\]

Making use of the result \cite[p.~141]{DLMF}:
\[\frac{\g(k-\sigma)}{\g(k)}=k^{-\sigma}\left\{1+\frac{\sigma(\sigma+1)}{2k}+O(k^{-2})\right\}\qquad (k\ra+\infty),\]
we therefore obtain\footnote{We use the shorthand notation $O(A, B, C)$ to represent $O(A)+O(B)+O(C)$.} for fixed $x$
\[{}_1\Psi_1(\rho,k;\rho,0;x)=\frac{xk^\rho \g(k)}{\g(\rho)}\left\{1-\frac{\rho(1-\rho)}{2k}+O(k^{-2})+
\frac{xk^\rho\g(\rho)}{2\g(2\rho)}(1+O(k^{-1}))\right.\]
\[\left.\hspace{8cm}+\frac{x^2k^{2\rho} \g(\rho)}{6\g(3\rho)}(1+O(k^{-1}))+O\left(\frac{k^{3\rho}}{\g(4\rho)}\right)\right\}\]
\bee\label{e41}
=\frac{xk^\rho \g(k)}{\g(\rho)}\left\{1+
\frac{xk^\rho\g(\rho)}{2\g(2\rho)}-\frac{\rho(1-\rho)}{2k}+\frac{x^2k^{2\rho} \g(\rho)}{6\g(3\rho)}+O\left(k^{-2},\frac{k^{\rho-1}}{\g(2\rho)},\frac{k^{3\rho}}{\g(4\rho)}\right)\right\}
\ee
as $k\ra+\infty$ for fixed $x\neq 0$ and $\rho\in(-1,0)$.

It is found numerically that (\ref{e41}) remains a good approximation when $k\ra+\infty$ and $k|x|^{1/\rho}\ra+\infty$. In the case of negative $x$, this numerical observation can also be derived from a local large deviation limit theorem for partial sums of independent integer-valued random variables in the case of power tails of the common probability function of these variables, which is given in \cite[Theorem 3.7.1]{BB}. 
The leading term in the above approximation, namely
\bee\label{e42}
{}_1\Psi_1(\rho,k;\rho,0;x)\sim\frac{xk^\rho \g(k)}{\g(\rho)}\qquad (k\ra+\infty,\ \rho\in(-1,0)),
\ee
was given previously without proof 
by Panjer and Willmot \cite[Ex.~10.2.2, p.~342]{PW}. We remark that when $\rho=-\fs$, the refining term involving $k^\rho$ inside the curly braces of the expansion in (\ref{e41}) vanishes; similarly, the order terms involving $k^{\rho-1}/\g(2\rho)$ and $k^{3\rho}/\g(4\rho)$ also vanish to yield that
\bee\label{e42a}
{}_1\Psi_1(-\fs,k;-\fs,0;x)=-\frac{x\,\g(k)}{2\sqrt{\pi k}}\bl\{1+\bl(\frac{3}{8}-\frac{x^2}{4}\bl)\frac{1}{k}+O(k^{-2})\br\}\qquad(k\ra+\infty).
\ee
The leading term in (\ref{e42a}) can be seen to agree with that obtained from (\ref{e22a}) and (\ref{e22b})
upon use of the large-$k$ asymptotic forms of the Bessel functions given in \cite[Section 10.41(i)]{DLMF}.
\vspace{0.2cm}

\noindent
\vspace{0.2cm}

\noindent
4.2.2\ {\it Asymptotics for integer $k\ra+\infty$ and $x\ra-\infty$}
\vspace{0.1cm}

\noindent
In this section, we replace negative $x$ by $-x$, where $x>0$. We derive an integral representation  for ${}_1\Psi_1(\rho, k; \rho,0;-x)$ analogous to (\ref{e400}) when $\rho=-\sigma$, $0<\sigma<1$ and $k$ is a positive integer.

Now, for $k=0, 1, 2, \ldots\,$, we have
 \[\frac{\g(k-\sigma n)}{k! \g(-\sigma n)}=\frac{1}{2\pi i}\int^{(1+)} \frac{t^{k-\sigma n-1}}{(t-1)^{k+1}}dt,\]
where integration is taken around a positively orientated closed path surrounding the point $t=1$, but not enclosing the origin. This result is easily established by evaluating the residue of the integrand at $t=1$. Then, from (\ref{e14}), it follows that
\begin{eqnarray}
{}_1\Psi_1(-\sigma,k;-\sigma,0;-x)&=&\frac{k!}{2\pi i} \int^{(1+)} \frac{t^{k-1}}{(t-1)^{k+1}}\sum_{n=0}^\infty
\frac{(-xt^{-\sigma})^n}{n!}\,dt \nonumber\\
&=&\frac{k!}{2\pi i}\int^{(1+)} \frac{t^{k-1} e^{-xt^{-\sigma}}}{(t-1)^{k+1}}\,dt,\label{e450}
\end{eqnarray}
where the integration path is as above.

Suppose that $x=O(k)$, and set $k:=ux$, where $u>0$ is fixed. Then (\ref{e450})  can be rewritten as
\[{}_1\Psi_1(-\sigma,k;-\sigma,0;-x)=\frac{k!}{2\pi i}\int^{(1+)} \frac{e^{k\psi(t)}}{t(t-1)}\,dt,\qquad\psi(t)=\log\,\frac{t}{t-1}-ut^{-\sigma}.\]
The phase function $\psi(t)$ has a saddle point $t_s$ on the positive real $t$-axis given by the unique
root of the equation
\bee\label{e451}
t_s^{-\sigma} (t_s-1)=\frac{1}{\sigma u},
\ee
which lies in $(1,\infty)$.\footnote{With $t=1/y$, $a=\sigma u>0$, (\ref{e451}) can be expressed in the form $y^{1+\rho}=a(1-y)$ when $\rho\in (-1,0)$. Then, by \cite[Proposition 5]{BKKK}, the unique root $y_s$ on $(0,1)$ admits the following closed-form representation in terms of the `reduced' Wright function $\phi$ for $\rho\in(-\frac{1}{2},0)$:
\[y_s=\bl\{-\log \bl(\int_0^\infty \frac{e^{-ay}}{y(1-y)}\, \phi(-(1+\rho),0;\frac{y-1}{(ay)^{\rho+1}})dy\br)\br\}^{\!1/(\rho+1)}.\]
When $\rho=-\frac{1}{2}$, it is straightforward to verify that $y_s=(2a/(1+\sqrt{1+4a^2}))^2$, which is consistent with the above formula. We do not yet have a proof when $\rho\in (-1,-\frac{1}{2})$, although we have been able to confirm this representation numerically.}
The path of steepest descent through $t_s$ crosses the real axis in a perpendicular direction and terminates at the origin in a similar manner to that shown in Fig.~1(a).
Following the same procedure as described in Section 4.1.1 (we omit the details), we therefore obtain that
\bee\label{e452}
{}_1\Psi_1(-\sigma, k; -\sigma,0; -x)\sim 
\g(k)\sqrt{\frac{k}{2\pi(1-\sigma)}} \left(t_{s0}+\frac{\sigma}{1-\sigma}\right)^{\!\!-1/2} \!\exp\, [k\psi(t_{s})]\,\{1+O(k^{-1})\}
\ee
as integer $k\ra+\infty$ and $x=O(k)$, where 
\bee\label{e453}
\psi(t_{s0})=-\log\,(1-\tau_s)-\frac{\tau_s}{\sigma(1-\tau_s)},\qquad \tau_s=1/t_{s}.
\ee
The value of the saddle $t_{s}\equiv t_{s}(u)>1$ is determined from (\ref{e451}). 

It is seen that the leading term in (\ref{e452}) is of the same form as that in (\ref{e45}) when $\rho>0$ (with $\rho$ replaced by $-\sigma$), albeit
with a different determination of the saddle point $t_s$.
In addition, it is not difficult to establish from (\ref{e453}) that $\psi(t_{s0})<0$, with the result that the factor $\exp\,[k\psi(t_{s0})]$ is exponentially small as $k\ra+\infty$.

\vspace{0.6cm}

\begin{center}
{\bf 5. \ Asymptotics of $\phi(\rho,0;x)$ for $x\ra\pm\infty$}
\end{center}
\setcounter{section}{5}
\setcounter{equation}{0}
\renewcommand{\theequation}{\arabic{section}.\arabic{equation}}
In \cite{W1, W2}, Wright studied the asymptotic expansion of what we have called in this paper the `reduced' Wright function (also known as a generalized Bessel function) 
${}_0\Psi_1(-\!\!\!-;\rho,\delta;z)\equiv \phi(\rho,\delta;z)$
for $|z|\ra\infty$ when the parameter $\rho>0$ and $\rho\in(-1,0)$.
We recall from (\ref{e13}) that the particular case of this function corresponding to $\delta=0$ is given by
\bee\label{e51}
{}_0\Psi_1(-\!\!\!-;\rho,0;x)=\phi(\rho,0;x)=\sum_{n=1}^\infty \frac{x^n}{n! \g(\rho n)}~,
\ee
where the variable $x$ is real. 

In this section, we present a summary of the known asymptotics of this function as $x\ra\pm\infty$ when $\rho\in (-1,0)\cup (0,\infty)$; see also \cite[Appendix A.2]{VPY2}. It is found that when $\rho>0$, the behavior of $\phi(\rho,0;x)$ is exponentially large as $x\ra+\infty$, but can be exponentially large, oscillatory or exponentially decaying as $x\ra-\infty$ according to the value of $\rho$. When $\rho\in(-1,0)$ and $x\ra+\infty$, the behavior depends sensitively on the value of the parameter $\rho$, whereas when $x\ra-\infty$ the behavior is exponentially decaying for all values of $\rho\in (-1,0)$.
\vspace{0.2cm}

\noindent 5.1\ {\it The expansion of $\phi(\rho,0;x)$ as $x\ra\pm\infty$ when $\rho>0$}
\vspace{0.1cm}

\noindent
The expansion of $\phi(\rho,0;x)$ when $\rho>0$ is obtained from \cite{W1} in the form
\bee\label{e501}
\phi(\rho,0;x)\sim\left\{\begin{array}{ll} E(x) & (x\ra+\infty)\\
\\
E(xe^{\pi i})+E(xe^{-\pi i}) & (x\ra-\infty),\end{array}\right.
\ee
where
\bee\label{e501a}
E(x):=\bl(\frac{\rho}{2\pi\kappa}\br)^\fr (hx)^{1/(2\kappa)} \exp\, [\kappa (hx)^{1/\kappa}] \sum_{j=0}^\infty c_j(\rho)(\kappa(hx)^{1/\kappa})^{-j}
\ee
and, from (\ref{e12}) with $\alpha=0$, $\beta=\rho$, we have $\kappa=1+\rho$ and $h=\rho^{-\rho}$; see also \cite[\S 2.3]{PK}. The coefficients $c_j(\rho)$ (with $c_0(\rho)=1$)
are discussed in  \cite[Appendix A.3]{VPY2}. If we set
\bee\label{e502}
c_j(\rho)=\frac{(2\rho+1)(\rho+2)}{2^{3j}\,\rho^j\,j!}\,\wp_j(\rho) \quad(j\geq 1),
\ee
where $\wp_j(\rho)$ is a polynomial in $\rho$ of degree $2(j-1)$, the first few polynomials are given by
\[\wp_1(\rho)=-\frac{1}{3},\qquad \wp_2(\rho)=\frac{1}{9}(2-19\rho+2\rho^2),\] \[\wp_3(\rho)=\frac{1}{135}(556+1628\rho-9093\rho^2+1628\rho^3+556\rho^4),\]
\[\wp_4(\rho)=-\frac{1}{405}(4568-226668\rho-465702\rho^2+2013479\rho^3-465702\rho^4-226668\rho^5+4568\rho^6).\]
The polynomials $\wp_j(\rho)$ are related to the so-called (V.M.) Zolotarev polynomials; see \cite[pp.~456--457]{VPY1} and \cite[Appendix A.3]{VPY2}.

When $\rho>1$, the expansion of $\phi(\rho,0;x)$ also includes additional exponential expansions of the same type as $E(x)$, but with the argument of $x$ rotated by multiples of $2\pi$. However, as these additional expansions are subdominant compared to $E(x)$ on the positive $x$-axis, we can still say that the dominant expansion is as stated in (\ref{e501}).

From (\ref{e501}), it can be seen that $\phi(\rho,0;x)$ is exponentially large as $x\ra+\infty$ for fixed $\rho>0$.
As $x\ra-\infty$, the growth of $\phi(\rho,0;x)$ is controlled by the exponential factors $\exp [\kappa(hx)^{1/\kappa} e^{\pm\pi i/\kappa}]$ appearing in $E(xe^{\pm\pi i})$, which depend on the sign of $\cos\,(\pi/\kappa)=\cos\,(\pi/(1+\rho))$. Thus, when $0<\rho<1$ we have $\cos\,(\pi/\kappa)<0$ and $\phi(\rho,0;x)$ decays exponentially as $x\ra-\infty$. When $\rho=1$, $\cos\,(\pi/\kappa)=0$ and $\phi(\rho,0;x)$ is oscillatory with an amplitude growing like $x^{1/4}$. Finally, when $\rho>1$ we have 
$\cos\,(\pi/\kappa)>0$ and $\phi(\rho,0;x)$ is consequently exponentially large as $x\ra-\infty$.
\vspace{0.2cm}

\noindent 5.2\ {\it The expansion of $\phi(\rho,0;x)$ as $x\ra\pm\infty$ when $\rho\in (-1,0)$}
\vspace{0.1cm}

\noindent
To discuss the case $\rho\in (-1,0)$ we set $\rho=-\sigma$, where $0<\sigma<1$. 
Then use of the reflection formula for the gamma function shows that
\bee\label{e52}
\phi(-\sigma,0;x)=-\frac{1}{\pi}\sum_{n=1}^\infty\frac{x^n}{n!}\,\g(1+\sigma n)\,\sin \pi\sigma n
=\frac{1}{2\pi i}\{F(xe^{-\pi i\sigma})-F(xe^{\pi i\sigma})\},
\ee
where
\bee\label{e520}
F(z)\equiv{}_1\Psi_0(\sigma,1;-\!\!\!-;z):=\sum_{n=0}^\infty \frac{z^n}{n!}\,\g(1+\sigma n)\qquad (0<\sigma<1).
\ee
The asymptotic expansion of $\phi(-\sigma,0;x)$ as $x\ra\pm\infty$ can then be obtained from that of $F(z)$ for $|z|\ra\infty$, which we describe below following \cite{Br}, \cite[\S 2.3]{PK}; see also \cite[\S 3]{PLMJ}.

From (\ref{e12}), we have the parameters associated with the function $F(z)$ in (\ref{e520}) given by
$\kappa=1-\sigma$ and $h=\sigma^\sigma$.
We introduce the formal algebraic and exponential expansions $H_{1,0}(z)$ and $E_{1.0}(z)$, respectively, by
\bee\label{e52a}
H_{1,0}(z):=\frac{1}{\sigma} \sum_{k=0}^\infty \frac{(-)^k}{k!}\,\g\left(\frac{k+1}{\sigma}\right)\,z^{-(k+1)/\sigma}
\ee
and
\bee\label{e52b}
E_{1,0}(z):=A_0' Z^{1/2}\,e^Z \sum_{j=0}^\infty c_j(-\sigma) Z^{-j},\qquad Z=\kappa(hz)^{1/\kappa},
\ee
where $A_0'=(2\pi\sigma)^{1/2}/\kappa$ and the coefficients $c_j(-\sigma)$ are defined in (\ref{e502}).
Then, from \cite{Br}, \cite{PLMJ}, \cite[\S 2.3]{PK}, we have the asymptotic expansion
\bee\label{e54}
F(z)\sim\left\{\begin{array}{lll} E_{1,0}(z)+H_{1,0}(ze^{\mp\pi i})& \mbox{in} & |\arg\,z|<\pi\kappa\\
\\
H_{1,0}(ze^{\mp\pi i}) & \mbox{in} & |\arg (-z)|<\pi(1-\kappa)\end{array}\right.
\ee
as $|z|\ra\infty$, where the upper or lower sign is chosen according as $z$ is situated in the upper or lower half-plane, respectively. 

In the sector $|\arg\,z|<\fs\pi\kappa$, $E_{1,0}(z)$ is exponentially large and makes the dominant contribut\-ion\footnote{The subdominant algebraic expansion $H_{1,0}(ze^{\mp\pi i})$ undergoes a Stokes phenomenon on $\arg\,z=0$.} to $F(z)$; see Fig.~2. The expansion $E_{1,0}(z)$, which is oscillatory on the rays $\arg\,z=\pm\fs\pi\kappa$ (the anti-Stokes lines), is still present beyond the sector $|\arg\,z|<\fs\pi\kappa$ where it becomes subdominant in the sectors $\fs\pi\kappa<|\arg\,z|<\pi\kappa$. The rays $\arg\,z=\pm\pi\kappa$, where $E_{1,0}(z)$ is {\it maximally}
subdominant with respect to $H_{1,0}(ze^{\mp\pi i})$, are called {\it Stokes lines}. As these rays are crossed (in the sense of increasing $|\arg\,z|$) the exponential expansion switches off according to the now familiar error-function smoothing law (see \cite[p.~67]{DLMF} for an overview of this topic) to leave the algebraic expansion in the sector $|\arg (-z)|<\pi(1-\kappa)=\pi\sigma$.
\begin{figure}[t]
\centering
\begin{picture}(200,200)(-100,-100)
\put(-85,0){\line(1,0){170}}
\put(0,-85){\line(0,1){170}}
\put(0,0){\line(3,2){80}}
\put(0,0){\line(3,-2){80}}
\put(0,0){\line(2,5){38}}
\put(0,0){\line(2,-5){38}}
\put(45,95){Stokes line}
\put(45,85){$\theta=\pi\kappa$}
\put(45,-85){Stokes line} 
\put(45,-95){$\theta=-\pi\kappa$}
\put(85,50){$\theta=\fs\pi\kappa$}
\put(85,-57){$\theta=-\fs\pi\kappa$}
\put(30,52){ES + A}
\put(30,-52){ES + A}
\put(70,20){EL + A}
\put(70,-20){EL + A}
\put(-40,30){A}
\put(-40,-30){A}
\end{picture}
\caption{\footnotesize{The exponentially large (EL), exponentially small (ES) and algebraic (A) sectors associated with $F(z)$ in the complex $z$-plane when $0<\kappa<1$ and $\theta=\arg\,z$. On the rays $\arg\,z=\pm\fs\pi\kappa$ the algebraic and exponential (oscillatory) expansions are of comparable significance.}}
\end{figure}
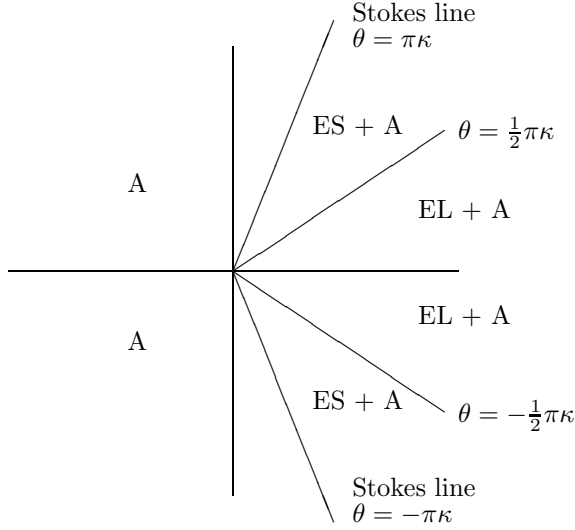
The expansion on the rays $\arg\,z=\pm\pi\kappa$ has been discussed in \cite[\S 5]{PLMJ}; see also below. 
\vspace{0.5cm}

\noindent 5.2.1\ {\it The expansion of $\phi(-\sigma,0;x)$ as $x\ra +\infty$}
\vspace{0.2cm}

\noindent
From (\ref{e52}) and (\ref{e52a}), the algebraic expansion associated with $\phi(-\sigma,0;x)$ is given by
\[H(x)\equiv \frac{1}{2\pi i}\{H_{1,0}(xe^{-\pi i\sigma}e^{\pi i})-H_{1,0}(xe^{\pi i\sigma} e^{-\pi i})\}
=\frac{1}{\pi\sigma}\sum_{k=0}^\infty\frac{x^{-(k+1)/\sigma}}{k!}\,\g\left(\frac{k+1}{\sigma}\right)\,\sin \frac{\pi(k+1)}{\sigma}\]
\bee\label{e54r}
=\frac{1}{\sigma}\sum_{k=0}^\infty\frac{x^{-(k+1)/\sigma}}{k! \g\left(1-\frac{k+1}{\sigma}\right)}.
\ee
From (\ref{e52}) and (\ref{e52b}), some routine algebra shows that the exponential expansion associated with $\phi(-\sigma,0;x)$ is 
\[E(x)\equiv\frac{1}{2\pi i}\{E_{1,0}(xe^{-\pi i\sigma})-E_{1,0}(xe^{\pi i\sigma})\}\]
\bee\label{e55}
=\frac{A_0'}{\pi}\,X^{1/2} e^{X \cos (\pi\sigma/\kappa)} \sum_{j=0}^\infty c_j(-\sigma) X^{-j} \cos\left\{X \sin \frac{\pi\sigma}{\kappa}+\frac{\pi}{\kappa}\left(\frac{1}{2}-\sigma j\right)\right\},
\ee
where 
\bee\label{e55a}
X:=\kappa(hx)^{1/\kappa}
\ee
and we recall that $\kappa=1-\sigma$, $h=\sigma^\sigma$.

Let us denote the points $z=xe^{\pm\pi i\sigma}$ appearing as the arguments of the function $F(z)$ in (\ref{e52}) by $P_\pm$. Then, when $0<\sigma<\f{1}{3}$, $P_\pm$ lie in the exponentially large sector of Fig.~2, and consequently as $x\ra+\infty$ we have
\bee\label{e56a}
\phi(-\sigma,0;x)\sim E(x)+H(x)\qquad (0<\sigma<\f{1}{3}).
\ee
Since $\cos\,(\pi\sigma/\kappa)>0$ for this range of $\sigma$, the asymptotic series $E(x)$ in (\ref{e55}) is exponentially large as $x\ra+\infty$ and is the dominant component of the expansion.
When $\sigma=\f{1}{3}$, $P_\pm$ lie on the anti-Stokes lines $\arg\,z=\pm\fs\pi\kappa$. On these rays, 
$\cos\,(\pi\sigma/\kappa)=0$ and the algebraic expansion in (\ref{e54r}) satisfies
$H(x)\equiv 0$, so that
\begin{eqnarray}
\phi(-\f{1}{3},0;x)&\sim&\frac{A_0'}{\pi}\,X^{1/2} \sum_{j=0}^\infty c_j(-\f{1}{3})X^{-j} \cos \left\{X+\frac{3\pi}{4}-\frac{\pi j}{2}\right\}\nonumber\\
&=&\frac{x^{3/4}}{3^{1/4}\sqrt{\pi}} \sum_{j=0}^\infty c_j(-\f{1}{3}) \bl(\frac{2x^{3/2}}{3\surd 3}\br)^{\!\!-j} \cos \bl\{\frac{2x^{3/2}}{3\surd 3}+\frac{3\pi}{4}-\frac{\pi j}{2}\br\}\label{e56b}
\end{eqnarray}
as $x\ra+\infty$. This case is related to the Airy Ai function (see (\ref{e61b})) and, from \cite[p.~198]{DLMF}, the coefficients $c_j(-\f{1}{3})$ can be expressed in closed form as  \[c_j(-\f{1}{3})=\frac{\g(3j+\fs)}{54^j j! \g(j+\fs)}.\]

When $\f{1}{3}<\sigma<\fs$, $P_\pm$ lie in the exponentially small and algebraic sectors of Fig.~2 and
\bee\label{e56c}
\phi(-\sigma,0;x)\sim E(x)+H(x)\qquad (\f{1}{3}<\sigma<\fs).
\ee
Since $\cos\,(\pi\sigma/\kappa)<0$ for this range of $\sigma$, $E(x)$ is exponentially small as $x\ra+\infty$ and is the subdominant component of the expansion.

When $\sigma=\fs$, $P_\pm$ lie on the Stokes lines $\arg\,z=\pm\fs\pi$. In this case, we have the exact result:
\bee\label{e56d}
\phi(-\fs,0;-x)=\frac{x}{2\sqrt{\pi}}\,e^{-x^2/4}\qquad (x\in {\bf R}).
\ee
We note that $\cos\,(\pi\sigma/\kappa)=-1$, $H(x)\equiv 0$ and, from (\ref{e502}), the coefficients $c_j(-\sigma)=0$ ($j\geq 1$) when $\sigma=\fs$.
On the Stokes lines $\arg\,z=\pm\pi\kappa$, the exponential expansion of $F(z)$ is equal to $\fs E_{1,0}(z)$ together with the addition of another exponentially small expansion (see below) which cancels in the linear combination (\ref{e52}). This yields the result 
\bee\label{e56e}
E(x)=\frac{A_0'}{2\pi} X^{1/2} e^{-X}\equiv\phi(-\fs,0;-x),
\ee
since $X=x^2/4$ when $\sigma=\fs$. The fact that in this case, the `reduced' Wright function coincides with the leading term of the corresponding exponentially small expansion is consistent with the well-known result that the saddle-point approximation for the probability density function of an inverse Gaussian distribution is exact, coinciding with this density (compare to \cite[Theorem~5.6]{BJ97}).

Finally, when $\fs<\sigma<1$, $P_\pm$ lie in the algebraic sector and accordingly
\bee\label{e56f}
\phi(-\sigma,0;x) \sim H(x)\qquad (\fs<\sigma<1)
\ee
as $x\ra+\infty$. The results in (\ref{e56a})--(\ref{e56f}) then describe the five different asymptotic forms of $\phi(-\sigma,0;x)$ as $x\ra+\infty$, when $\sigma=-\rho\in (0,1)$, and agree with those given by Wright in \cite{W2}.
\vspace{0.5cm}

\noindent 5.2.2\ {\it The expansion of $\phi(-\sigma,0;x)$ as $x\ra -\infty$}
\vspace{0.2cm}

\noindent
To deal with the expansion for $x\ra-\infty$, we replace $x$ by $-x$ to find from (\ref{e52}) that
\bee\label{e57}
\phi(-\sigma,0,-x)=\frac{1}{2\pi i}\{F(xe^{\pi i\kappa})-F(xe^{-\pi i\kappa})\}.
\ee
The arguments $xe^{\pm\pi i\kappa}$ are situated on the Stokes lines $\arg\,z=\pm\pi\kappa$, respectively, for all values of $\sigma$ satisfying $0<\sigma<1$. It is readily seen that 
\[H_{1,0}(xe^{\pi i\kappa-\pi i})-H_{1,0}(xe^{-\pi i\kappa+\pi i})\equiv 0,\]
so that the algebraic expansions present in the combination (\ref{e57}) cancel. It has been shown in \cite[Eq.~(4.24)]{PLMJ} that when the algebraic expansion $H_{1,0}(ze^{\mp\pi i})$ is optimally truncated at, or near, its least term, then the exponential expansion associated with $F(z)$ on the Stokes lines $\arg\,z=\pm\pi\kappa$ assumes the form
\[\frac{1}{2} E_{1,0}(xe^{\pm\pi i\kappa})-\sqrt{\frac{2}{\pi}}\,e^{-X}\sum_{j=0}^\infty (-)^j {\hat B}_j X^{-j}.\]
Here $X$ is defined in (\ref{e55a}) and the coefficients ${\hat B}_j$ (which we do not specify here) depend on the coefficients $c_j$. It can be shown (compare \cite[Section 5]{PLMJ}) that the series involving the coefficients ${\hat B}_j$ cancel in the combination in (\ref{e52}), and hence, we deduce that
\begin{eqnarray}
\phi(-\sigma,0;-x) &\sim& \frac{1}{4\pi i}\{E_{1,0}(xe^{\pi i\kappa})-E_{1,0}(xe^{-\pi i\kappa})\}\nonumber\\
&=&\frac{A_0'}{2\pi}\,X^{1/2} e^{-X} \sum_{j=0}^\infty c_j(-\sigma) X^{-j} \qquad (0<\sigma<1)\label{e58}
\end{eqnarray}
as $x\ra+\infty$. This expansion is seen to be exponentially small for large $x>0$ when $0<\sigma<1$ and was obtained recently in \cite[Eq.~(A.14)]{VPY2}. The leading term in (\ref{e58}) has been derived previously by Mikusi\'nski \cite{JM} who employed different methods. 

It is also relevant that a somewhat related discussion in a probabilistic context can be found in \cite[Sections 2.10--2.11]{Z}. For instance, the leading behavior on the right-hand side of (\ref{e56b}) can be derived with some effort from a combination of (\ref{e61b}) with \cite[p.~171]{Z}.
\vspace{0.6cm}

\begin{center}
{\bf 6. \ Some special representations and properties of $\phi(\rho,0;x)$}
\end{center}
\setcounter{section}{6}
\setcounter{equation}{0}
\renewcommand{\theequation}{\arabic{section}.\arabic{equation}}
In this final section, we list some integral representations and structural properties of the `reduced' Wright function (or generalized Bessel function) $\phi(\rho,0;x)$ that we consider to be of interest. Most of the results we present relate to the case when $\rho=-\sigma$, $0<\sigma<1$.
\vspace{0.2cm}

\noindent 6.1\ {\it Special representations}
\vspace{0.1cm}

\noindent
We first note the special values of $\sigma$ for which $\phi(-\sigma,0;\pm x)$ can be expressed in terms of a well-known special function. These are:
\bee\label{e61g}
\phi(1,0;x)=x^{1/2} I_1(2x^{1/2}),\qquad \phi(1,0;-x)=-x^{1/2} J_1(2x^{1/2})\quad(x\geq 0),
\ee
\bee\label{e61a}
\phi(-\fs,0,x)= -\frac{x}{2\sqrt{\pi}}\,e^{-x^2/4}\qquad(x\in{\bf R}),
\ee
\begin{eqnarray}
\phi(-\f{1}{3},0,-x)&=&3^{-1/3}x\, \mbox {Ai}(3^{-1/3}x)\nonumber\\
&=&\left\{\begin{array}{ll}-\dfrac{(-x)^{3/2}}{3\surd 3}\bl\{J_\frac{1}{3}\bl(\dfrac{2(-x)^{3/2}}{3\surd 3}\br)+J_{-\frac{1}{3}}\bl(\dfrac{2(-x)^{3/2}}{3\surd 3}\bl)\br\} & (x\leq 0)\\ \\
\dfrac{x^{3/2}}{3\pi}\,K_\frac{1}{3}\bl(\dfrac{2x^{3/2}}{3\surd 3}\br)
\quad (x\geq0),\end{array}\right.\label{e61b}
\end{eqnarray}
\bee\label{e61e}
\phi(-\f{2}{3},0,x)=-\frac{1}{2\sqrt{3\pi}}\,e^{2x^3/27}\,W_{-\frac{1}{2},\frac{1}{6}}\bl(\frac{4x^3}{27}\bl)\quad (x\geq0),
\ee
\bee\label{e61c}
\phi(-\f{2}{3},0,-x)=\sqrt{\frac{3}{\pi}}\,e^{-2x^3/27}\,W_{\frac{1}{2},\frac{1}{6}}\bl(\frac{4x^3}{27}\br)\quad(x\geq0),
\ee
where $J_\nu$, $I_\nu$ and $K_\nu$ are the usual Bessel functions, and Ai and $W_{\kappa,\mu}$ are the Airy and Whittaker functions, respectively. The first part of expression (\ref{e61g}) was noted in \cite[erratum]{VPY2}.
The expressions (\ref{e61b}), corresponding to $x\geq 0$, (\ref{e61c}) and (\ref{e22b})  and the integral representations  (\ref{e23}), (\ref{e24}) are seen to imply with some effort special cases of the structural property (\ref{e15}) for negative $\rho$. Similarly,
the expression (\ref{e61a}), when combined with the integral representation for the $K$ Bessel function in \cite[Eq.~(10.32.10)]{DLMF}, is also an example of (\ref{e15}). 
The representations (\ref{e61b}) (for positive $x$) and (\ref{e61c}) have been discussed more fully in \cite[Eqs.~(4.9), (4.3)]{VPY2}; the derivation of (\ref{e61b})--(\ref{e61c}) is given in Appendix C.

\vspace{0.2cm}

\noindent 6.2\ {\it Integral representations}
\vspace{0.1cm}

\noindent
We start from the following integral representation (see \cite[Eq.~(2)]{RS} or \cite[Section 1.4]{W1}):
\bee\label{e62}
\phi(-\sigma,0;x)=\frac{1}{2\pi i}\int_C e^{\tau+x\tau^\sigma}d\tau \qquad (0<\sigma<1),
\ee
where $C$ is a path in the complex $\tau$-plane that starts at infinity in the third quadrant, passes round the origin in the positive sense and returns to infinity in the second quadrant. The representation (\ref{e62}) holds for complex $x$ satisfying $-\pi\leq\arg\,x\leq\pi$.

If the path $C$ is collapsed onto each side of the branch cut along the negative real axis together with a small circular path round the origin of vanishingly small radius, we obtain the following representation given in \cite{RS}:
\bee\label{e63}
\phi(-\sigma,0;\pm x)=\mp\frac{x^{-1/\sigma}}{\pi}\int_0^\infty \exp\,[-x^{-1/\sigma}t\pm t^\sigma \cos \pi\sigma]\,\sin (t^\sigma \sin \pi\sigma)\,dt
\ee
for $x>0$ and $0<\sigma<1$.
 
Another integral representation that can be obtained from \cite[Eq.~(2)]{JM} is
\bee\label{e64}
\phi(-\sigma,0;-x)=\frac{\sigma x^{1/(1-\sigma)}}{\pi(1-\sigma)} \int_0^\pi W_\sigma(\varphi) \exp\,[-x^{1/(1-\sigma)}W_\sigma(\varphi)]\,d\varphi\qquad(x>0;\ 0<\sigma<1),
\ee
where
\[W_\sigma(\varphi):=\bl(\frac{\sin \sigma\varphi}{\sin \varphi}\br)^{\sigma/(1-\sigma)}\,\frac{\sin (1-\sigma)\varphi}{\sin \varphi}.\]
An equivalent expression in the context of positive stable probability density functions can be derived from  \cite[Eq.~(2.5.10)]{Z}. This is explained in \cite[Eq.~(3.17)]{VPY1}, where a slightly different trigonometric function than $W_\sigma(\varphi)$ and integration range are used.

A simple deduction that can be made from this last representation concerns the positivity of $\phi(-\sigma,0;-x)$ on the positive real axis. Since $W_\sigma(\varphi)>0$ for $0<\varphi<\pi$, it follows trivially from (\ref{e64}) that
\bee\label{e65}
\phi(-\sigma,0;-x)>0 \qquad (x>0;\ 0<\sigma<1);
\ee
see \cite[p.~192]{JM}.
It is clear that the positivity stipulated by (\ref{e65}) provides an analytical justification for employing this `reduced' Wright function in probability theory\footnote{A combination of \cite[Eqs.~(3.9), (3.11), (3.15)--(3.16), Remark 2]{VPY1} and \cite[Eqs.~(2.3.2)--(2.3.3)]{Z} implies with some effort that $x^{-1} \phi(-\sigma,0;-x)>0$ for $\sigma\in[\frac{1}{2},1)$ and $x\in{\bf R}$.}. See also \cite[Remark 3]{VPY1} for an important relationship to stochastic modelling.
 
A representation of $\phi(-\sigma,0;-x)$ involving an integral of a similar function is given in the following theorem, which resembles \cite[Eq.~(2.10.11)]{Z} where an analogous formula for the Mittag-Leffler function is given.
\begin{theorem}$\!\!\!.$\ \  
For $x>0$ and $0<\sigma<1$, we have
\bee\label{e66}
\phi(-\fs\sigma,0;-x)=\frac{x}{2\sqrt{\pi}}\int_0^\infty s^{-3/2}e^{-x^2/(4s)}\,\phi(-\sigma,0;-s)\, ds.
\ee
\end{theorem}
\noindent Proof.\ \ In the integral representation (\ref{e62}), we first observe that it is always possible to choose the path $C$ in the $\tau$-plane such that $\Re (\tau^\sigma)>0$ when $0<\sigma<1$. 
If we substitute (\ref{e62}) for $\phi(-\sigma,0;-s)$ into the right-hand side of (\ref{e66}),  
we obtain that
\[\int_0^\infty s^{-3/2}e^{-x^2/(4s)}\,\phi(-\sigma,0;-s)\, ds=
\frac{1}{2\pi i}\int_C e^\tau \bl(\int_0^\infty e^{-s\tau^{\sigma}-x^2/(4s)} s^{-3/2}ds\br)\,d\tau\]
\[=\frac{2\sqrt{\pi}}{x}\,\frac{1}{2\pi i}\int_Ce^{\tau-x\tau^{\sigma/2}}d\tau
=\frac{2\sqrt{\pi}}{x}\,\phi(-\fs\sigma,0;-x)\]
upon reversal of the order of integration and use of the result \cite[Eq.~(10.32.10)]{DLMF}
\[\int_0^\infty e^{-s\tau^\sigma -x^2/(4s)} s^{-3/2}ds=\frac{2\sqrt{\pi}}{x}\,e^{-xs^{\sigma/2}},\]
which is valid when $\Re (\tau^\sigma)>0$.\hfill $\Box$
\vspace{0.3cm}

It is relevant that Theorem 2 constitutes an important special case of the subsequent Theorem 3. In addition, it implies new integral representations for two specific values of the `reduced' Wright function; see (\ref{e67}) and (\ref{e67c}).
Thus, if we put $\sigma=\fs$ in (\ref{e66}) and use (\ref{e61a}), we obtain that
\[\phi(-\f{1}{4},0;-x)=\frac{x}{2\sqrt{\pi}} \int_0^\infty s^{-3/2} e^{-x^2/(4s)}\,\phi(-\fs,0;-s)\,ds\]
\[=\frac{x}{4\pi}\int_0^\infty s^{-1/2} e^{-s^2/4-x^2/(4s)}\,ds.\]
With the change of variable $s\ra x^{2/3}\tau^2$, this yields the integral representation
\bee\label{e67}
\phi(-\f{1}{4},0;-x)=\frac{x^{4/3}}{2\pi}\int_0^\infty \exp\,[\f{1}{4}x^{4/3} (\tau^4+\tau^{-2})]\,d\tau.
\ee
The validity of (\ref{e67}) can also be derived with some effort from \cite[Eq.~(3.4.26)]{Z}.
Similarly, if we put $\sigma=\f{1}{3}$ in (\ref{e66}) and use (\ref{e61b}), we obtain that
\begin{eqnarray}
\phi(-\f{1}{6},0;-x)&=&\frac{3^{-1/3}x}{2\sqrt{\pi}} \int_0^\infty s^{-1/2}\,e^{-x^2/(4s)} \mbox{Ai}\,(3^{-1/3}s)\,ds\nonumber\\
&=&\frac{x}{6\pi^{3/2}} \int_0^\infty e^{-x^2/(4s)} K_\frac{1}{3}\bl(\frac{2s^{3/2}}{3\surd 3}\br)\,ds.\label{e67c}
\end{eqnarray}

Next, we have the Laplace transforms \cite[Eq.~(11)]{RS}
\bee\label{e67a}
\int_0^\infty e^{-zt} t^{\delta-1} \phi(-\sigma,\delta;-\alpha t^{-\sigma})\,dt=z^{-\delta} \exp\,(-\alpha z^\sigma)
\qquad (0<\sigma<1,\ \alpha>0)\ee
and, for $\rho>0$, $\alpha>0$
\bee\label{e67b}
\int_0^\infty e^{-zt} t^{\delta-1} \phi(\rho,\delta;\pm\alpha t^{\rho})\,dt=\left\{\begin{array}{ll}z^{-\delta} \exp\,(\pm\alpha z^{-\rho}) & (\delta\neq 0)\\ 
\\
\exp\,(\pm\alpha z^{-\rho})-1 & (\delta=0).\end{array}\right.
\ee
The result in (\ref{e67b}) is obtained by substitution of the series expansion of $\phi(-\sigma,\delta;\pm\alpha t^\rho)$ followed by term-by-term integration.
We also have the Mellin transform \cite[Eq.~(21)]{RS}
\bee\label{e68}
\int_0^\infty t^{\mu-1}\phi(-\sigma,0;-tx^{-\sigma})\,dt=\frac{x^{\sigma\mu}\g(\mu)}{\g(\sigma\mu)}\qquad (\mu>-1).
\ee
It is relevant that the evaluation (\ref{e68}) can also be derived with some effort from the moment representations of stable distributions (compare to \cite[Eq.~(2.11.11)]{Z}).

A combination of probabilistic and analytical arguments presented in \cite[Eqs.~(3.9), (3.11), (3.15)--(3.16), Remark 2]{VPY1} and \cite[Eqs.~(2.3.2)--(2.3.3), p.~178, Property 7b]{Z} implies that for $s\geq 0$
\bee\label{e68a}
\int_{-\infty}^\infty e^{st}t^{-1}\,\phi(-\sigma,0;-t)\,dt=\exp\,(s^{1/\sigma})
\ee
when $\sigma=[\f{1}{3},1)$; the above integral is divergent when $\sigma\in(0,\f{1}{3})$ on account of the exponential growth of $\phi(-\sigma,0,-t)$ given in (\ref{e56a}). This result has been verified numerically for
several $\sigma\in[\f{1}{3},1)$. In particular, if we put $\sigma=\f{1}{3}$ in (\ref{e68a}) and use (\ref{e61b}), then we obtain for $s\geq 0$ the evaluation
\bee\label{e68b}
3^{-1/3} \int_{-\infty}^\infty \mbox{Ai}\,(3^{-1/3}t)\,e^{st}dt=\exp\,(s^3).
\ee
A version of (\ref{e68b}) is given in \cite[p.~171]{Z}. It is also relevant that the functions $t^{-1} \phi(-\sigma,0;-t)$ with $\sigma\in(0,\fs)$ belong to the so-called class of `trans-stable' functions considered in \cite[p.~174]{Z}. It is also relevant to quote \cite[p.~259, Comment 2.35]{Z} and \cite[p.~255, Comment after Eq.~(2)]{OS}, which implicitly claims the validity of the Fourier-transform counterpart of (\ref{e68b}).

\vspace{0.2cm}

\noindent 6.3\ {\it New multiplication properties and the reflection principle}
\vspace{0.1cm}

\noindent
First, we establish three general multiplication properties that involve a `reduced' Wright function $\phi$ expressed in terms of an integral of a product of two $\phi$ functions. In view of (\ref{e61a}), the first one generalizes Theorem 2 above.

\begin{theorem}$\!\!\!.$\ \ 
Let $\sigma_1$ and $\sigma_2$ be such that $0<\sigma_1<1$ and $0<\sigma_2<1$. Then for $x>0$ we have
\bee\label{e69a}
\phi(-\sigma_1\sigma_2,0;-x)=\int_0^\infty \phi(-\sigma_1,0;-s)\,\phi(-\sigma_2,0;-xs^{-\sigma_2})\,\frac{ds}{s}~.
\ee
\end{theorem}

\noindent Proof.\ \ 
We substitute the integral representation (\ref{e62}) for $\phi(-\sigma_1,0;-s)$ on the right-hand side of (\ref{e69a}) to find that
\[I:=\int_0^\infty \phi(-\sigma_1,0;-s)\,\phi(-\sigma_2,0;-xs^{-\sigma_2})\,\frac{ds}{s}=
\frac{1}{2\pi i}\int_C e^\tau\bl(\int_0^\infty e^{-s\tau^{\sigma_1}}\phi(-\sigma_2,0;-xs^{-\sigma_2})\frac{ds}{s}\br) d\tau.\]
Making use of the Laplace transform result (\ref{e67a}) to evaluate the inner integral, we have
\[I=\frac{1}{2\pi i}\int_C e^{\tau-x\tau^{\sigma_1\sigma_2}} d\tau=\phi(-\sigma_1\sigma_2,0;-x)\]
by (\ref{e62}), which concludes the proof.\hfill$\Box$
\bigskip

An expression analogous to (\ref{e69a}) for the class of Mittag-Leffler functions is given in \cite[Theorem 2.10.3]{Z}.

\begin{theorem}$\!\!\!.$\ \
Let $\rho_1$ and $\rho_2$ be positive parameters such that either (i) $\rho_1\leq 1$, $\rho_2<1$ or (ii) $\rho_1<1$, $\rho_2\geq 1$ with $\rho_1 \rho_2<1$. Then for $x>0$ we have
\bee\label{e69b}
\phi(-\rho_1 \rho_2,0; -x)=\int_0^\infty \phi(\rho_1,0;-s)\,\phi(\rho_2,0; -xs^{\rho_2})\,\frac{ds}{s}.
\ee
\end{theorem}
Before giving the proof of Theorem 4, we first observe from the asymptotics of $\phi(\rho,0;z)$ in (\ref{e501}) and (\ref{e501a}) that the absolute value of the integrand in (\ref{e69b}) as $s\ra+\infty$ is controlled by the exponential factor
\bee\label{e62d}
\exp\,\bl[a_1 \cos\bl(\frac{\pi}{1+\rho_1}\br) s^{1/(1+\rho_1)}+a_2 \cos \bl(\frac{\pi}{1+\rho_2}\br) (xs^{\rho_2})^{1/(1+\rho_2)}\br],
\ee
where $a_j:=(1+\rho_j)^{-1}\rho_j^{-\rho_j/(1+\rho_j)}$ ($j=1, 2$). It is clear that, when $\rho_1\leq 1$, $\rho_2<1$, or $\rho_1<1$, $\rho_2\geq 1$ (with $\rho_1\rho_2<1$),
this factor decreases exponentially as $s\ra+\infty$, so that the integral in (\ref{e69b}) converges
absolutely. However, when $\rho_1>1$, $\rho_2>1$, or $\rho_1>1$, $\rho_2<1$ (with $\rho_1\rho_2<1$), the exponential factor (\ref{e62d}) increases without limit as $s\ra+\infty$, with the result that the integral in (\ref{e69b}) is not defined in these cases.
 
When $\rho_1=\rho_2=1$, the exponential factor (\ref{e62d}) equals 1 and the (improper) integral on the right-hand side of (\ref{e69b}) does not exist. This can be understood by use of the second relation in (\ref{e61g}). The integral in (\ref{e69b}) can then be expressed in terms of
\[\int_0^\infty t J_1(2t) J_1(2x^{1/2} t)\,dt.\] 
But in view of \cite[Eq.~(10.22.56)]{DLMF} this integral does not converge.
\bigskip

\noindent Proof.\ \ From the representation analogous to (\ref{e62}), which can be found in \cite[\S 1.4]{W1}, one obtains that for $\rho>0$
\bee\label{e69c}
\phi(\rho,0;-x)=\frac{1}{2\pi i} \int_C e^{\tau-x\tau^{-\rho}}\,d\tau,
\ee
where $C$ is the loop surrounding the negative real $\tau$-axis as before. The right-hand side of (\ref{e69b}) can therefore be written as
\[I_1=\frac{1}{2\pi i}\int_C e^\tau \bl(\int_0^\infty e^{-s\tau^{-\rho_1}}\,\phi(\rho_2,0; -xs^{\rho_2})\,\frac{ds}{s}\br) d\tau.\]

The inner integral may be evaluated by use of the Laplace transform in (\ref{e67b}) to yield that
\[I_1=\frac{1}{2\pi i} \int_C  e^\tau \{e^{-x\tau^{\rho_1\rho_2}}-1\}\,d\tau=\phi(-\rho_1\rho_2,0;-x)\qquad (\rho_1\rho_2<1)\]
by (\ref{e62}), where the $-1$ in the integrand makes no contribution since $\int_C e^\tau d\tau=0$.
The condition $\rho_1\rho_2<1$, which is clearly satisfied when (i) or (ii) is fulfilled, applies provided the integral in (\ref{e69b}) converges. This yields the conditions stated in the theorem.\hfill $\Box$ 

\begin{theorem}$\!\!\!.$\ \ 
Let $\sigma$ and $\rho$ be such that $0<\sigma<1$ and $\rho>0$. Then for $x>0$ we have
\bee\label{e69z}
\phi(\rho\sigma,0;\pm x)=\int_0^\infty \phi(-\sigma,0;-s)\,\phi(\rho,0;\pm xs^{\rho})\,\frac{ds}{s}~.
\ee
\end{theorem}

\noindent Proof.\ \ 
It follows the same procedure employed in establishing Theorem 3 and uses (\ref{e62}), the Laplace transform result (\ref{e67b}) and, finally, the representation (\ref{e69c}). \hfill $\Box$ 
\bigskip

\noindent{\bf Corollary 1.}\ \ {\it For $\sigma\in(0,1)$ and $x>0$, the following reflection principle holds:}
\bee\label{e70e}
\left.\begin{array}{rcl}
\phi(\sigma,0;x)&=&\displaystyle{\sqrt{x}\int_0^\infty \phi(-\sigma,0; -s)\,s^{-1/2} I_1(2\sqrt{xs})\,ds}\\
\\
\phi(-\sigma,0;-x)&=&-\displaystyle{\sqrt{x}\int_0^\infty \phi(\sigma,0; -s)\,s^{-1/2} J_1(2\sqrt{xs})\,ds}.
\end{array}\right\}
\ee
\vspace{0.2cm}

\noindent Proof.\ \ 
If we put $\rho=1$ in (\ref{e69z}), we find that for $x>0$
\bee\label{e70d}
\phi(\sigma,0;\pm x)=\int_0^\infty \phi(-\sigma,0;-s)\,\phi(1,0;\pm xs)\,\frac{ds}{s}\qquad (0<\sigma<1).
\ee
The result (\ref{e70e}) follows immediately by taking the upper signs in (\ref{e70d}) and using (\ref{e69b}) with $\rho_1=\sigma\in (0,1)$, $\rho_2=1$ combined with the Bessel function representations in (\ref{e61g}).\hfill $\Box$
\bigskip

Some other integrals involving products of `reduced' Wright functions are cited in \cite{RS}. We remark that (\ref{e69a}) can also be obtained with some effort from \cite[Eq.~(3.3.1)]{Z} which states\footnote{The sign ${\stackrel{d}{=}}$ means that the distributions of random variables coincide.} that
$Y_1(\sigma_1) Y_2(\sigma_2){\stackrel{d}{=}} Y(\sigma_1 \sigma_2)$,
where $Y_1(\sigma_1)$ and $Y_2(\sigma_2)$ are two independent positive stable random variables with corresponding indices of stability $\sigma_1$ and $\sigma_2$, and $Y(\sigma_1 \sigma_2)$ is the positive stable random variable with index of stability $\sigma_1\sigma_2$ such that the corresponding p.d.f.'s of these random variables are given by $x^{-1}\phi(-\sigma,0;-x^{-\sigma})$ (compare to \cite[Eq.~(3.10)]{VPY1}).

A combination of the representations (\ref{e61g})--(\ref{e61b}) and (\ref{e61c}) with $\sigma=\fs$, $\f{1}{3}$ and $\f{2}{3}$ then leads to the following integral representations in terms of products of well-known special functions:
\bee\label{e70a}
\phi(\fs,0;\pm x)=\pm \frac{1}{2}\sqrt{\frac{x}{\pi}} \int_0^\infty s^{1/2}e^{-s^2/4} {\cal C}_1(2\sqrt{sx})\,ds,
\ee
\bee\label{e70b}
\phi(\f{1}{3},0,\pm x)=\pm \frac{\sqrt{x}}{3\pi}\int_0^\infty s\,K_\frac{1}{3}\bl(\frac{2s^{3/2}}{3\surd 3}\br)\,{\cal C}_1(2\sqrt{xs})\,ds
\ee
and
\bee\label{e70c}
\phi(\f{2}{3},0,\pm x)=\pm \sqrt{\frac{3x}{\pi}}\int_0^\infty s^{-1/2} e^{-2s^3/27} W_{\frac{1}{2},\frac{1}{6}}\bl(\frac{4s^3}{27}\br)\,{\cal C}_1(2\sqrt{xs})\,ds
\ee
for $x>0$,
where ${\cal C}_1$ denotes the Bessel functions $I_1$ (upper signs) and $J_1$ (lower signs). 

We regard Corollary 1 and its illustrations given by (\ref{e70a})--(\ref{e70c}) to be of interest in their own right. However, it is yet unclear whether these formulas might have a computational value. For instance, a substitution of the well-known series representations for $I_1$ and $J_1$ into (\ref{e70a}) shows that
\[\phi(\fs,0; \pm x)=\frac{1}{2\sqrt{\pi}}\sum_{n=1}^\infty \frac{\g(\fs n+\fs) (\pm 2x)^n}{\g(n)\,n!},\]
which is comparable to the original series representation of this function given by (\ref{e51}).

Our final theorem concerns an asymptotic connection between the Wright function ${}_1\Psi_1$ and its reduced counterpart $\phi$.
\begin{theorem}$\!\!\!.$\ \ 
For $\rho>0$ and $x>0$, we have as $k\ra+\infty$
\bee\label{e610}
{}_1\Psi_1(\rho,k; \rho,0;x) \sim \g(k)\,\phi(\rho,0;xk^\rho).
\ee
The same result holds when $\rho\in(-1,0)$ and $x\in{\bf R}$ as $k\ra+\infty$.
\end{theorem}

\noindent Proof.\ \ Consider first the case with $\rho>0$ and $x>0$. With $\kappa=1+\rho$ and $h=\rho^\rho$, the leading behavior of ${}_1\Psi_1(\rho,k; \rho,0;x)$ from (\ref{e48}) and (\ref{e48a}) can be expressed in the following form:
\[{}_1\Psi_1(\rho,k; \rho,0;x)\sim \g(k)\,\bl(\frac{\rho}{2\pi\kappa}\br)^{\!\!\fr} (hxk^\rho)^{1/2\kappa} \exp\,[\kappa(hxk^\rho)^{1/\kappa}]\]
as $k\ra+\infty$ with $k/x\ra+\infty$. But from (\ref{e501}), the leading behavior of $\phi(\rho,0;xk^\rho)$ is
\[\bl(\frac{\rho}{2\pi\kappa}\br)^{\!\!\fr} (hxk^\rho)^{1/2\kappa} \exp\,[\kappa(hxk^\rho)^{1/\kappa}]\]
as $xk^\rho\ra+\infty$, which implies (\ref{e610}).

When $\rho\in (-1,0)$, we have from (\ref{e42}) with $\rho=-\sigma$, $0<\sigma<1$ and $x\in {\bf R}$ that
\[{}_1\Psi_1(-\sigma,k; -\sigma,0;x) \sim \frac{xk^{-\sigma}\g(k)}{\g(-\sigma)}\qquad (k\ra+\infty).\]
From the series representation for $\phi(-\sigma,0;xk^{-\sigma})$ in (\ref{e51}), it is readily seen that as $k\ra+\infty$, the argument $xk^{-\sigma}\ra 0$ for fixed $x$, and the first term yields the leading behavior.  Hence, we have 
\[\phi(-\sigma,0;xk^{-\sigma})\sim \frac{xk^{-\sigma}}{\g(-\sigma)}\qquad (k\ra+\infty)\]
and the result follows.\hfill $\Box$
\vspace{0.2cm}

To conclude this section, we note that an asymptotic relationship of the same type as (\ref{e610}) in the case where $\rho\in(-1,-\fs]$, $k\ra+\infty$ and $x\ra-\infty$ can be derived with some effort from a probabilistic local limit theorem stipulated by \cite[Eq.~(1.8)]{KN}. (Note that the case $\rho\in(-\fs,0)$ should be excluded
in order to have the finiteness of the so-called first pseudomoment; compare \cite[Section 3.2]{V94}.)
\vspace{0.6cm}

\begin{center}
{\bf Appendix A: \ Proof of the structural representation (\ref{e15})}
\end{center}
\setcounter{section}{1}
\setcounter{equation}{0}
\renewcommand{\theequation}{\Alph{section}.\arabic{equation}}
In this appendix, we establish the structural representation given in (\ref{e15}) that relates ${}_1\Psi_1$ to the `reduced' Wright function $\phi$. This can be separated into two cases according as $\rho$ is positive or negative as follows: 
\bee\label{b1}\left.
\begin{array}{rl}	
{}_{1}\Psi_{1}(\rho, k; \rho, 0; x)&=
\displaystyle{\int_{0}^{\infty}}~e^{-\tau}\tau^{k-1} \phi(\rho, 0; x \tau^{\rho})\,d\tau \quad (\rho>0,\ k\geq 0)\\
\\	
{}_{1}\Psi_{1}(\rho, k; \rho, 0; -x)&=
\displaystyle{\int_{0}^{\infty}}~e^{-\tau}\tau^{k-1}\ \phi(\rho, 0; -x \tau^{\rho})\,d\tau \quad (\rho\in (-1,0),\ k=0, 1, 2, \ldots).
\end{array}\right\}
\ee

The case when $\rho>0$ is straightforward. Substitution of the series expansion for $\phi(\rho, 0;x\tau^\rho)$ in (\ref{e51}), followed by reversal of the order of summation and integration, yields that the right-hand side of the first expression in (\ref{b1}) can be written in the following form:
\[\sum_{n=1}^\infty \frac{x^n}{n!\,\g(\rho n)} \int_0^\infty e^{-\tau} \tau^{k+\rho n-1}d\tau
=\sum_{n=1}^\infty \frac{\g(k+\rho n)}{\g(\rho n)}\,\frac{x^n}{n!}\]
valid for $k\geq 0$. By (\ref{e14}), this is seen to equal the left-hand side of (\ref{b1}).
We note that the case $k=0$ is trivial, since ${}_1\Psi_1(\rho,0;\rho,0;x)=e^x$ and, from (\ref{e67a}), the first integral in (\ref{b1}) also reduces to $e^x$.

When $\rho\in (-1,0)$ we put $\rho=-\sigma$ in the second expression in (\ref{b1}), with $0<\sigma<1$. We then have to show that
\bee\label{b3}
{}_1\Psi_1(-\sigma, k;-\sigma, 0; -x)=\int_0^\infty e^{-\tau} \tau^{k-1} \phi(-\sigma,0;-x\tau^{-\sigma})\,d\tau
\ee
for non-negative integer $k$. Now, for $k=0, 1, 2, \ldots\,$, we have
\begin{eqnarray}
{}_1\Psi_1(-\sigma,k;-\sigma,0;-x)&=&\sum_{n=1}^\infty \frac{\g(k-\sigma n)}{\g(-\sigma n)}\,\frac{(-x)^n}{n!}\nonumber\\
&=&(-)^k\sum_{n=1}^\infty \frac{(-x)^n}{n!}\,\sigma n(\sigma n-1) \ldots (\sigma n-k+1),\label{b4}
\end{eqnarray}
with the product being empty (equal to 1) when $k=0$.
But from (\ref{e67a}) we have the Laplace transform for $z>0$
\[\int_0^\infty e^{-z\tau} \tau^{-1}\,\phi(-\sigma,0; -x\tau^{-\sigma})\,d\tau=\exp\,(-xz^\sigma).\]
Compare to \cite[p.~452]{VPY1} where a connection to a positive stable density function and its Laplace transform is made; see also a comment above Theorem 6.
Differentiation of this result $k$ times then yields that
\begin{eqnarray}
\int_0^\infty e^{-\tau} \tau^{k-1} \,\phi(-\sigma,0;-x\tau^{-\sigma})\,d\tau&=&(-)^k\,\frac{d^k}{dz^k}(e^{-xz^\sigma})|_{z=1}
=(-)^k \sum_{n=1}^\infty \frac{(-x)^n}{n!}\,\frac{d^k}{dz^k} (z^{\sigma n})|_{z=1}\nonumber\\
&=&(-)^k \sum_{n=1}^\infty \frac{(-x)^n}{n!}\,\sigma n(\sigma n-1)\ldots (\sigma n-k+1).\label{b5}
\end{eqnarray}

The evaluation in (\ref{b5}) is seen to equal the series representation (\ref{b4}) of the left-hand side of (\ref{b3}) when $k=0, 1, 2, \ldots\,$, thereby concluding the proof; compare to  \cite[Theorem 1a]{H}. \hfill $\Box$

\vspace{0.6cm}

\begin{center}
{\bf Appendix B: \ Proof of the integral representations in Theorem 1}
\end{center}
\setcounter{section}{2}
\setcounter{equation}{0}
\renewcommand{\theequation}{\Alph{section}.\arabic{equation}}
Denote the integral on the right-hand side of (\ref{e23}) by $I_3$. 
From the Mellin-Barnes integral representation for the modified Bessel function of the second kind (see, for example, \cite[Eq.~(10.32.13)]{DLMF}) we have
\[K_\frac{1}{3}(z)=\frac{(\fs x)^{1/3}}{4\pi i}\int_{c-\infty i}^{c+\infty i} \g(s) \g(s-\f{1}{3}) (\fs z)^{-2s}ds\qquad (c>\f{1}{3},\ |\arg\,z|<\fs\pi).\]
With $z=a/\sqrt{t}$, $a=2(x/3)^{3/2}$ ($x>0$), we then find that
\begin{eqnarray*}
I_3&=&\frac{x^{3/2}}{3\pi} \int_0^\infty e^{-t} t^{k-3/2} K_\frac{1}{3}(a/\sqrt{t})\,dt\\
&=&
\frac{x^2}{3^{3/2}\pi}\int_0^\infty e^{-t}t^{k-5/3}\left\{\frac{1}{4\pi i}\int_{c-\infty i}^{c+\infty i} \g(s) \g(s-\f{1}{3}) (x/3)^{-3s} t^sds\right\}dt\\
&=&\frac{x^2}{6\surd 3\,\pi}\cdot\frac{1}{2\pi i}\int_{c-\infty i}^{c+\infty i} \g(s) \g(s-\f{1}{3}) \g(k+s-\f{2}{3}) (x/3)^{-3s} ds
\end{eqnarray*}
upon reversal of the order of integration and evaluation of the integral over $t$ as a gamma function.
Making use of the triplication formula 
\bee\label{etf}
\g(s) \g(s-\f{1}{3}) \g(s-\f{2}{3})=\frac{2\pi}{3^{2s-5/2}} \g(3s-2)
\ee
for the gamma function (see, for example, \cite[Eq.~(5.5.6)]{DLMF}), we obtain that
\[I_3=\frac{3x^2}{2\pi i} \int_{c-\infty i}^{c+\infty i} \g(3s-2)\,
\frac{\g(k+s-\f{2}{3})}{\g(s-\f{2}{3})}\, x^{-3s}ds.\]

For $k=1, 2, \ldots\,$, the integrand has poles only at $s=\f{2}{3}-\f{1}{3}n$ ($n=1, 2, \ldots$) resulting from $\g(3s-2)$, since the quotient of gamma functions is regular. Displacement of the integration path to the left to coincide with the vertical line $\Re (s)=\f{5}{6}-\f{1}{3}N$, passing over the first $N-1$ poles of $\g(3s-2)$ with residue $\f{1}{3}(-)^n/ n!$, we find
\bee\label{e25}
I_3=\sum_{n=1}^{N-1} \frac{\g(k-\f{1}{3}n)}{\g(-\f{1}{3}n)}\,\frac{(-x)^n}{n!}+R_N,
\ee
where, on putting $s=\f{5}{6}-\f{1}{3}N+\f{1}{3}it$ with $-\infty<t<\infty$, 
\[R_N=\frac{x^{N-1/2}}{2\pi} \int_{-\infty}^\infty \g(-N\!+\!\fs\!+\!it)\,{\cal G}_N(t)x^{-it}\,dt, \qquad {\cal G}_N(t)=\prod_{m=0}^{k-1} (s-\f{2}{3}+m).\]

Now, $|s-\f{2}{3}+m|\leq\f{1}{3}(N^2+t^2)^{1/2}\leq \f{1}{3}N(1+t^2)^{1/2}$ provided $N>3m>1$, so that
\[|{\cal G}_N(t)|\leq (N/3)^k (1+t^2)^{k/2}\qquad (N>3k).\]
By repeated application of the well-known formula $z\g(z)=\g(z+1)$, we see that
\[|\g(-N\!+\!\fs\!+\!it)|=\frac{|\g(\fs+it)|}{|(N\!-\!\fs\!-\!it)(N\!-\!\f{3}{2}\!-\!it)\dots (\fs\!-\!it)|}\leq \frac{\sqrt{\pi}\,|\g(\fs+it)|}{\g(N+\fs)}.\]
With these bounds we then have
\[|R_N|\leq \frac{x^{N-1/2} (N/3)^k}{2\sqrt{\pi}\,\g(N+\fs)} \int_{-\infty}^\infty |\g(\fs+it)| (1+t^2)^{k/2}dt,\]
where the integral is convergent, since $|\g(\fs+it)|=(\pi/\cosh \pi t)^{1/2}$.
It therefore follows that
$R_N=O(N^k x^{N-1/2})/\g(N+\fs)$ and hence, that $R_N\ra 0$ as $N\ra\infty$ for fixed $x>0$. We can then put $N=\infty$ in (\ref{e25}), which is seen, by (\ref{e14}), to be the series expansion of ${}_1\Psi_1(-\f{1}{3},k;-\f{1}{3},0;-x)$. This establishes (\ref{e23}) for positive integer $k$.

A similar procedure can be employed to establish (\ref{e24}). We have the Mellin-Barnes integral representation
given in \cite[Eq.~(13.16.12)]{DLMF}
\[e^{-z/2} W_{\fr,\frac{1}{6}}(z)=\frac{1}{2\pi i}\int_{c-\infty i}^{c+\infty i} \frac{\g(s+\f{1}{3}) \g(s+\f{2}{3})}{\g(s+\fs)}\,z^{-s}ds\qquad (c>0,\ |\arg\,z|<\fs\pi).\]
With $z=a/t^2$, $a=4x^3/27$ ($x>0$), the integral on the right-hand side of (\ref{e24}) may be expressed as
\begin{eqnarray*}
I_4&=&\sqrt{\frac{3}{\pi}}\int_0^\infty e^{-t}t^{k-1}\exp(-\fs a/t^2) W_{\frac{1}{2}, \frac{1}{6}}(a/t^2)dt\\
&=&\sqrt{\frac{3}{\pi}}\,\frac{1}{2\pi i} \int_{c-\infty i}^{c+\infty i} \frac{\g(s+\f{1}{3}) \g(s+\f{2}{3})}{\g(s+\fs)}\,\g(k+2s) \left(\frac{4x^3}{27}\right)^{\!-s} ds\\
&=& \frac{3}{2\pi i} \int_{c-\infty i}^{c+\infty i}\g(3s)\,\frac{\g(k+2s)}{\g(2s)}\,x^{-3s}ds
\end{eqnarray*}
upon use of the duplication and triplication formulas for the gamma function; compare to (\ref{etf}).

Provided $k=1, 2, \ldots \,$, the only poles of the above integrand result from the factor $\g(3s)$ situated at $s=-\f{1}{3}n$ ($n=1, 2, \ldots$).
Displacement of the integration path to the left to coincide with the vertical line $\Re(s)=-\f{1}{3}N+\f{1}{6}$ and evaluation of the residues then yields that
\bee\label{e26}
I_4=\sum_{n=1}^{N-1}\frac{\g(k-\f{2}{3}n)}{\g(-\f{2}{3}n)}\,\frac{(-x)^n}{n!}+R_N,
\ee
where, by similar arguments to those above, $R_N=O(N^k x^{N-1/2})/\g(N+\fs)\ra 0$ as $N\ra\infty$. Hence, we can put $N=\infty$ in (\ref{e26}), which then equals the series expansion of ${}_1\Psi_1(-\f{2}{3},k;-\f{2}{3},0;-x)$. This establishes (\ref{e24}) for positive integer $k$. \hfill $\Box$
\vspace{0.6cm}

\begin{center}
{\bf Appendix C: \ Proof of the special function representations in Section 6.1}
\end{center}
\setcounter{section}{3}
\setcounter{equation}{0}
\renewcommand{\theequation}{\Alph{section}.\arabic{equation}}
From (\ref{e51}) and (\ref{e52}) with $x\geq 0$, we have 
\[\phi(-\f{1}{3},0;\pm x)=-\frac{1}{\pi}\sum_{n=1}^\infty \frac{(\pm x)^n}{n!}\,\g(1+\f{1}{3}n)\,\sin (\f{1}{3}\pi n).\]
Replacing $n$ by $3m+j$, $j=0, 1, 2$, we obtain
\[\phi(-\f{1}{3},0;\pm x)=\mp \frac{x\surd 3}{2\pi}\bl\{\sum_{m=0}^\infty \frac{\g(m+\f{4}{3})}{(3m+1)!}\,(\mp x)^{3m} \pm x \sum_{m=0}^\infty \frac{\g(m+\f{5}{3})}{(3m+2)!}\,(\mp x)^{3m}\br\},\]
since the sum corresponding to $j=0$ vanishes.
Use of the triplication formula (see (\ref{etf})) for the gamma function and introduction of the new variable ${\cal X}\equiv x^3/27$ leads to
\begin{eqnarray}
\phi(-\f{1}{3},0;\pm x)&=&\mp\frac{x}{3}\bl\{\sum_{m=0}^\infty \frac{(\mp {\cal X})^m}{m! \g(m+\f{2}{3})}\pm {\cal X}^{1/3}\sum_{m=0}^\infty \frac{(\mp {\cal X})^m}{m! \g(m+\f{4}{3})}\br\}\nonumber\\
&=&\mp \frac{x}{3}\bl\{{\cal X}^{1/6} {\cal C}_{-\frac{1}{3}}(2\sqrt{{\cal X}})+{\cal X}^{1/6} {\cal C}_{\frac{1}{3}}(2\sqrt{{\cal X}})\br\},\label{c5}
\end{eqnarray}
where ${\cal C}_\nu(z)$ denotes the Bessel functions $J_\nu(z)$ (upper signs) and $I_\nu(z)$ (lower signs).

Taking the upper and lower signs in (\ref{c5}) separately, we therefore find that
\bee\label{c6a}
\phi(-\f{1}{3},0; x)=-\frac{x^{3/2}}{3\surd 3}\bl\{J_{-\frac{1}{3}}\bl(\frac{2x^{3/2}}{3\surd 3}\br)+
J_{\frac{1}{3}}\bl(\frac{2x^{3/2}}{3\surd 3}\br)\br\}\qquad (x\geq 0)
\ee
and
\begin{eqnarray}
\phi(-\f{1}{3},0;-x)&=&\frac{x^{3/2}}{3\surd 3}\bl\{I_{-\frac{1}{3}}\bl(\frac{2x^{3/2}}{3\surd 3}\br)-
I_{\frac{1}{3}}\bl(\frac{2x^{3/2}}{3\surd 3}\br)\br\}\nonumber\\
&=&\frac{x^{3/2}}{3\pi} K_{\frac{1}{3}}\bl(\frac{2x^{3/2}}{3\surd 3}\br)\qquad (x\geq 0).\label{c6b}
\end{eqnarray}
The representation of (\ref{c6a}) and (\ref{c6b}) in terms of the Ai Airy function follows upon use of \cite[Eqs.~(9.6.1), (9.6.6)]{DLMF}. This completes the proof of (\ref{e61b}).

The proof of (\ref{e61e}) and (\ref{e61c})  follows the same procedure. From (\ref{e52}), we have
\[\phi(-\f{2}{3},0;\pm x)=-\frac{1}{\pi}\sum_{n=1}^\infty \frac{(\pm x)^n}{n!}\,\g(1+\f{2}{3}n)\,\sin (\f{2}{3}\pi n),\]
which becomes, on replacing $n$ by $3m+j$, $j=0, 1, 2$,
\[\phi(-\f{2}{3},0;\pm x)=\mp \frac{x\surd 3}{2\pi}\bl\{\sum_{m=0}^\infty \frac{\g(2m+\f{5}{3})}{(3m+1)!}\,(\pm x)^{3m} \mp x \sum_{m=0}^\infty \frac{\g(2m+\f{7}{3})}{(3m+2)!}\,(\pm x)^{3m}\br\}.\]
Use of the duplication and triplication formulas for the gamma function and introduction of the new variable ${\cal X}\equiv 4x^3/27$ then leads to
\begin{eqnarray}
\phi(-\f{2}{3},0;\pm x)&=&\mp\frac{2^{2/3}x}{3\sqrt{\pi}}\bl\{\sum_{m=0}^\infty\frac{\g(m+\f{5}{6})}{\g(m+\f{2}{3})}\,(\pm {\cal X})^m
\mp z^{1/3} \sum_{m=0}^\infty \frac{\g(m+\f{7}{6})}{\g(m+\f{4}{3})}\,(\pm {\cal X})^m\br\}\nonumber\\
&=&\mp\frac{2^{2/3}x}{3\sqrt{\pi}}\bl\{\frac{\g(\f{5}{6})}{\g(\f{2}{3})}\,{}_1F_1(\f{5}{6};\f{2}{3};\pm{\cal X})\mp {\cal X}^{1/3}\,\frac{\g(\f{7}{6})}{\g(\f{4}{3})}\,{}_1F_1(\f{7}{6};\f{4}{3};\pm{\cal X})\br\}\nonumber\\
&=&\mp\frac{2^{2/3}x}{6\sqrt{3\pi}}\bl\{\frac{\g(\f{1}{3})}{\g(\f{7}{6})}\,{}_1F_1(\f{5}{6};\f{2}{3};\pm{\cal X})\pm {\cal X}^{1/3}\,\frac{\g(-\f{1}{3})}{\g(\f{5}{6})}\,{}_1F_1(\f{7}{6};\f{4}{3};\pm{\cal X})\br\}.\label{c1}
\end{eqnarray}

From \cite[Eqs. (13.2.42), (13.14.5)]{DLMF}, we have the following combination of confluent hypergeometric functions expressed in terms of the Whittaker function $W_{\kappa,\mu}(z)$ given by
\[\frac{\g(1-b)}{\g(a-b+1)}\,{}_1F_1(a;b;z)+z^{1-b}\,\frac{\g(b-1)}{\g(a)}\,{}_1F_1(a-b+1;2-b;z)\hspace{2cm}\]
\bee\label{c2}
\hspace{7cm}=e^{z/2} z^{b/2}W_{\frac{1}{2}b-a, \frac{1}{2}b-\frac{1}{2}}(z)
\ee
for $z\in C$.
Inserting the values $a=\f{5}{6}$, $b=\f{2}{3}$ in (\ref{c2}), we find upon taking the upper signs in 
(\ref{c1}) that
\begin{eqnarray}
\phi(-\f{2}{3}, 0;x)&=&-\frac{2^{2/3} x}{6\sqrt{3\pi}}\, e^{\fr {\cal X}} {\cal X}^{-1/3}\,W_{-\frac{1}{2}, \frac{1}{6}}({\cal X})\nonumber\\
&=&-\frac{1}{2\sqrt{3\pi}}\, e^{2x^3/27}\,W_{-\frac{1}{2}, \frac{1}{6}}\bl(\frac{4x^3}{27}\br)\quad (x\geq 0),\label{c3}
\end{eqnarray}
where we have made use of the fact that $W_{\kappa,-\mu}(z)=W_{\kappa, \mu}(z)$. This completes the proof of (\ref{e61e}).

Taking the lower signs in (\ref{c1}), we have
\[\phi(-\f{2}{3}, 0;-x)=\frac{2^{2/3}x}{\sqrt{3\pi}}\,e^{-{\cal X}}\bl\{\frac{\g(\f{1}{3})}{\g(\f{1}{6})}\,{}_1F_1(-\f{1}{6};\f{2}{3};{\cal X})+ {\cal X}^{1/3}\,\frac{\g(-\f{1}{3})}{\g(-\f{1}{6})}\,{}_1F_1(\f{1}{6};\f{4}{3};{\cal X})\br\},\]
upon use of the well-known Kummer transformation ${}_1F_1(a;b;z)=e^z\,{}_1F_1(b-a;b;-z)$. Insertion of the values $a=-\f{1}{6}$, $b=\f{2}{3}$ in (\ref{c1}) then yields that
\begin{eqnarray}
\phi(-\f{2}{3}, 0;-x)&=&\frac{2^{2/3} x}{\sqrt{3\pi}}\, e^{-{\cal X}/2} {\cal X}^{-1/3}\,W_{\frac{1}{2}, \frac{1}{6}}({\cal X})\nonumber\\
&=&\sqrt{\frac{3}{\pi}}\, e^{-2x^3/27}\,W_{\frac{1}{2}, \frac{1}{6}}\bl(\frac{4x^3}{27}\br)\quad (x\geq 0).\label{c4}
\end{eqnarray}
This completes the proof of (\ref{e61c}).
\hfill $\Box$

Although we do not use this result here, it is worth recording the connection between the Whittaker functions $W_{\pm\frac{1}{2}, \frac{1}{6}}(z)$. From \cite[Eqs.~(13.15.8), (13.18.0)]{DLMF}, we obtain the relation
\[W_{-\frac{1}{2}, \frac{1}{6}}(z)=6W_{\frac{1}{2}, \frac{1}{6}}(z)-\frac{6z}{\sqrt{\pi}}\,K_{\frac{1}{3}}(\fs z)\qquad (z\in C).\]

\vspace{0.6cm}

\noindent{\bf Acknowledgment:}\ \ V.V. wishes to thank the Fields Institute, University of Toronto and York University for their hospitality. He also acknowledges financial support from the Ohio University International Travel Fund for a visit to the first author in the UK.

\end{document}